\documentclass[11pt,epsfig,amsfonts]{amsart}

\newtheorem{df}{Definition}[section]
\newtheorem{lm}{Lemma}[section]
\newtheorem{thm}{Theorem}[section]

\newtheorem{cor}{Corollary}[section]
\numberwithin{equation}{section}

\usepackage{xcolor}
\usepackage[colorlinks,
            linkcolor=red,
            anchorcolor=blue,
            citecolor=green
            ]{hyperref}

\usepackage{epsfig}

\subjclass[2020]{Primary:  35B41, 37L30, 37L55; Secondary: 35B40.}

\keywords{McKean-Vlasov equation,  uniform measure attractor, complete solution, unbounded thin domain.
\\ This work was supported by  NSFC  (12371178 and 11971394) and Natural Science Foundation of Sichuan province under grant (2023NSFSC1342). All correspondences should be addressed to Li Ran.}

\author{ Tianhao Zeng }
\address[Tianhao Zeng ]
{ School of Mathematics \\
 Southwest Jiaotong University, Chengdu 610031,  China}
\email[T.~Zeng]{zengtianhao123@my.swjtu.edu.cn}

\author{ Ran Li}
\address[Ran Li]
{School of Mathematics \\
Southwest Jiaotong University, Chengdu 610031,  China}
\email[Ran.~Li]{liranzym@my.swjtu.edu.cn}

\author{Dingshi Li}
\address[Dingshi Li]
{School of Mathematics \\
Southwest Jiaotong University, Chengdu 610031,  China}
\email[D.~Li]{lidingshi@swjtu.edu.cn}

\begin{document}

\begin{abstract}
This article addresses the issue of uniform measure attractors for non-autonomous McKean-Vlasov stochastic reaction-diffusion equations defined on unbounded thin domains.
Initially, the concept of uniform measure attractors is recalled, and thereafter, the existence and uniqueness of such attractors are demonstrated.
Uniform tail estimates are employed to establish the asymptotic compactness of the processes,
thereby overcoming the non-compactness issue inherent in the usual Sobolev embedding on unbounded thin domains.
 Finally, we demonstrate that  the upper semi-continuity of uniform measure attractors defined on $(n+1)$-dimensional unbounded thin domains collapsing into the space $\mathbb R^n$.
\end{abstract}

\title[uniform measure  attractors of McKean-Vlasov stochastic reaction-diffusion equations]{
uniform measure  attractors of McKean-Vlasov stochastic reaction-diffusion equations
 on unbounded thin domain}
\maketitle

\section{introduction}
In this paper, our objective is to investigate the limiting behaviour of uniform measure attractors for the following non-autonomous, almost periodic stochastic reaction-diffusion equation driven by nonlinear noise defined on  $\mathcal O_{\varepsilon}$.
\begin{equation}\label{a1}
\left\{\begin{array}{l}
d \hat{u}^{\varepsilon}(t)-\Delta\hat{u}^{\varepsilon}(t)dt+\lambda \hat{u}^{\varepsilon}(t)dt+f\left( x, \hat{u}^{\varepsilon}(t,x),\mathcal L_{\hat{u}^{\varepsilon}(t)}\right) d t=g(t, x^*) d t \\
\quad+\sum\limits_{k=1}^{\infty}\left(\sigma_{k}(x)+\kappa(x^*)\varpi_k\left(\hat{u}^{\varepsilon}(t,x),\mathcal L_{\hat{u}^{\varepsilon}(t)}\right)\right) d W_k(t), \quad x\in\mathcal O_{\varepsilon},\ t>\tau, \\
\frac{\partial \hat{u}^{\varepsilon}}{\partial \nu_{\varepsilon}}=0,\quad x\in\partial\mathcal O_{\varepsilon},
\end{array}\right.
\end{equation}
with initial data
 \begin{equation}\label{a5}
\hat{u}^{\varepsilon}(\tau,x)=\hat{\xi}^{\varepsilon}(x),\quad x\in\mathcal O_{\varepsilon},
 \end{equation}
where $\lambda>0$ denotes a positive constant, $\mathcal L_{\hat{u}^{\varepsilon}(t)}$ represents the probability distribution of $\hat{u}^{\varepsilon}(t)$, $f$  is a nonlinear function possessing an arbitrary growth rate, and $g$ signifies an almost-periodic external term in $t$. Furthermore,
$\nu_{\varepsilon}$ denotes the unit outward normal vector to the boundary $\partial\mathcal O_{\varepsilon}$,
$\kappa\in L^2(\mathbb R^n)\cap L^{\infty}(\mathbb R^n)$, $\sigma_k$ is defined on $\widetilde{{\mathcal O}}$,
 $\varpi_k$ is a nonlinear diffusion term, and $(W_k)_{k\in\mathbb N}$ is a sequence of independent two-sided real-valued
 standard Wiener processes on a complete filtered probability space $(\Omega,\mathcal F,\{\mathcal F_t\}_{t\in\mathbb R},P)$ which satisfies the usual condition.

The unbounded thin domain, denoted as $\mathcal O_{\varepsilon}$, is rigorously defined as
$$
\mathcal{O}_{\varepsilon}=\left\{x=\left(x^*, x_{n+1}\right) \mid x^*=\left(x_1, \ldots, x_n\right) \in \mathbb{R}^n \text { and } 0<x_{n+1}<\varepsilon \rho\left(x^*\right)\right\},
$$
where $\rho \in C^2\left(\mathbb{R}^n,(0,+\infty)\right)$, and $0<\varepsilon \leq 1$.  Consequently, there exist constants $\rho_1$ and $\rho_2$ such that
\begin{equation}\label{a4}
\rho_1 \leq \rho\left(x^*\right) \leq \rho_2, \quad \forall x^* \in \mathbb{R}^n .
\end{equation}
Let $\mathcal{O}=\mathbb{R}^n \times(0,1)$ and $\widetilde{\mathcal{O}}=\mathbb{R}^n \times\left[0, \rho_2\right]$, where $\widetilde{\mathcal{O}}=\mathbb{R}^n \times\left[0, \rho_2\right]$ contains $\mathcal{O}_{\varepsilon}$ for $0<\varepsilon \leq 1$.

It is noteworthy that the unbounded thin domain $\mathcal O_{\varepsilon}$ undergoes a collapse, converging to the space $\mathbb R^n$ as $\varepsilon\rightarrow 0$. We will
investigate the existence and uniqueness and the limit of uniform measure attractors of \eqref{a1}-\eqref{a5} as $\varepsilon\rightarrow0$.

 For $\varepsilon=0$, the limit equation of \eqref{a1} reduces to the follow system defined on $\mathbb R^n$
\begin{align}\label{a2}
 \mathrm{d} u^0(t)-\frac{1}{\rho} &\sum_{i=1}^n\left(\rho u_{y_i}^0\right)_{y_i} \mathrm{~d} t+\lambda u^0(t)dt+f\left(\left(y^*, 0\right), u^0,\mathcal L_{u^0(t)}\right)dt=g\left(t,(y^*)\right) \mathrm{d} t \nonumber\\
&+\sum_{k=1}^{\infty}\left(\sigma_k\left(y^*, 0\right)+\kappa\left(y^*\right) \varpi_k\left(u^0(t),\mathcal L_{u^0(t)}\right)\right) \mathrm{d} W_k(t), \quad   t>\tau,
\end{align}
with initial condition
\begin{equation}\label{a3}
u^0\left(\tau, y^*\right)=\xi^0\left(y^*\right),
\end{equation}
where $y^*=(y_1,\cdots,y_n)\in\mathbb R^n$.

The McKean-Vlasov stochastic differential equations (MVSDEs) were initially examined in the seminal works of M. McKean  \cite{M67}
and V. Vlasov  \cite{V68}. These equations frequently emerge from the realm of interacting particle systems, as evidenced in the extensive literature\cite{BH,DV,FG}.
It is the distinctive quality of these differential equations that they depend not only on the states of the solutions, but also on the distributions of the solutions.
Consequently, the Markov operators pertinent to MVSDEs cease to constitute semigroups
are no longer semigroups (as illustrated in \cite{W2018}), thereby rendering the methodologies devised for addressing stochastic equations independent of distributions inapplicable to MVSDEs in a straightforward manner.

At present, a considerable number of academic publications are dedicated to the examination of solutions to MVSDEs.
For instance, the existence of solutions to MVSDEs has been examined in previous research  \cite{AD,FHSY,HDS}.  Additionally,  the existence and ergodicity of invariant measures for finite-dimensional MVSDEs have been investigated in other studies \cite{BSY,HW,Z23}. In a recent publication, Shi et al. \cite{SSLW} innovatively utilized the theory of pullback measure attractors to demonstrate the existence of invariant measures and periodic measures for infinite-dimensional MVSDEs.

For non-autonomous random dynamical
systems, there are typically two kinds of pathwise pullback random attractor which have drawn
much attention in last years: cocycle attractors introduced in \cite{W2012} and uniform random attractors
introduced in  \cite{CL}. For further details, please refer to the following sources:\cite{CK15,HUK,HK24,LCW,LI}.

 It bears noting that the notion of a measure attractor for autonomous dynamical systems in measure spaces was initially posited by Schmalfu$\ss$ in  \cite{S}.
For further details concerning the existence of measure attractors for autonomous stochastic equations, the reader is directed to the references \cite{MC,M,S2}.
A recent investigation by Li and Wang \cite{LW24} explored the limiting dynamical behaviour of pullback measure attractors in a system of semilinear parabolic stochastic equations with deterministic non-autonomous forcing in bounded thin domains. Notably, this investigation did not assume any particular constraints on the external forces involved, including translation-bounded. Subsequently, Li et al. initially explored uniform measure attractors in the context of stochastic Navier-Stokes equations,  as outlined in their work cited in \cite{LLZ}. In this paragraph, however, it is worth noting that in all of the articles mentioned above, the stochastic equations do not depend on the distributions of the solutions.

The dynamics of deterministic partial differential equations (PDEs) on thin domains were initially investigated by Hale and Raugel \cite{HR01,HR02,HR3}.
Subsequently, a considerable number of models have extended these initial findings, as demonstrated by the following references \cite{ACG,JKN,PR03,PR01}.
Recently, several results have been published concerning stochastic dynamical systems on thin domains. For example, see the references \cite{LWW17,LLWW18,LLWW19}
for studies of random attractors on bounded thin domains and the references \cite{SWLW19,SLLW20} for studies of random attractors on unbounded thin domains.

In order to illustrate the principal conclusions of this paper and for the sake of convenience, we shall use equations \eqref{c8}, which are  the equivalent of equations \eqref{a1}-\eqref{a5}.  we introduce the notation $\mathcal{P}\left(L^2\left(\mathcal O\right)\right)$ to represent the space of probability measures on $\left(L^2\left(\mathcal O\right), \mathcal{B}\left(L^2\left(\mathcal O\right)\right)\right.$, where $\mathcal{B}\left(L^2\left(\mathcal O\right)\right)$
 denotes the Borel
 $\sigma$-algebra associated with  $L^2\left(\mathcal O\right)$. Notably, the weak topology of  $\mathcal{P}\left(L^2\left(\mathcal O\right)\right)$ is metrizable, and the corresponding metric is denoted by $d_{\mathcal{P}\left(L^2\left(\mathcal O\right)\right)}$. Set

$$
\mathcal{P}_4\left(L^2\left(\mathcal O\right)\right)=\left\{\mu \in \mathcal{P}\left(L^2\left(\mathcal O\right)\right): \int_{L^2\left(\mathcal O\right)}\|\xi\|_{L^2\left(\mathcal O\right)}^4 d \mu(\xi)<\infty\right\}.
$$

Then $\left(\mathcal{P}_4\left(L^2\left(\mathcal O\right)\right), d_{\mathcal{P}\left(L^2\left(\mathcal O\right)\right)}\right)$ is a metric space. Given $r>0$, denote by

$$
B_{\mathcal{P}_4\left(L^2\left(\mathcal O\right)\right)}(r)=\left\{\mu \in \mathcal{P}_4\left(L^2\left(\mathcal O\right)\right): \int_{L^2\left(\mathcal O\right)}\|\xi\|_{L^2\left(\mathcal O\right)}^4 d \mu(\xi) \leq r^4\right\}.
$$

Given $\tau \leq t$ and $\mu \in \mathcal{P}_4\left(L^2\left(\mathcal O\right)\right)$, let $P_*^{g,\varepsilon}(t,\tau) \mu$ be the law of the solution of \eqref{c8} with initial law $\mu$ at initial time $\tau$. If $\phi: L^2\left(\mathcal O\right) \rightarrow \mathbb{R}$ is a bounded Borel function, then we write

$$
P^{g,\varepsilon}(\tau, t) \phi\left(\xi^{\varepsilon}\right)=\mathbb{E}\left(\phi\left(u^{\varepsilon}\left(t, \tau, \xi^{\varepsilon}\right)\right)\right), \quad \forall \xi^{\varepsilon} \in L^2\left(\mathcal O\right),
$$
where $u^{\varepsilon}\left(t, \tau, \xi^{\varepsilon}\right)$ is the solution of \eqref{c8} with initial value $\xi^{\varepsilon}$ at initial time $\tau$. Note that for the McKean-Vlasov stochastic equations like \eqref{c8}, $P_*^{g,\varepsilon}(t,\tau)$ is not the dual of $P^{g,\varepsilon}(\tau, t)$ (see \cite{W2018}) in the sense that

\begin{align}\label{a6}
\int_{L^2\left(\mathcal O\right)} P^{g,\varepsilon}(t,\tau) \phi(\xi) d \mu(\xi) \neq \int_{L^2\left(\mathcal O\right)} \phi(\xi) d P^{g, \varepsilon}_*(\xi),
\end{align}
where $\mu \in \mathcal{P}_4\left(L^2\left(\mathcal O\right)\right)$ and $\phi: L^2\left(\mathcal O\right) \rightarrow \mathbb{R}$ is a bounded Borel function.

To prove the existence of uniform measure attractors for \eqref{c8}, we first need to show $\left\{P^{g,\varepsilon}_*(t,\tau)\right\}_{\tau \leq t}$ is a
jointly continuous process, i.e., it is continuous in  $\mathcal P_4(L^2(\mathcal O))\times \mathcal H(g_0)$.
In general, if a stochastic equation does not depend on distributions of the solutions, then the continuity of $\left\{P^{g,\varepsilon}_*(t,\tau)\right\}_{\tau \leq t}$
on $\left(\mathcal{P}_4\left(L^2\left(\mathcal O\right)\right), d_{\mathcal{P}\left(L^2\left(\mathcal O\right)\right)}\right)$ follows from the Feller property of $\left\{P^{g,\varepsilon}(t,\tau)\right\}_{\tau \leq t}$ and the duality relation between $P^{g,\varepsilon}_*(t,\tau)$ and $P^{g,\varepsilon}(t,\tau)$. However, this method does not apply to the McKean-Vlasov stochastic reaction-diffusion equation \eqref{a1}-\eqref{a5}, because $P^{g,\varepsilon}_*(t,\tau)$ is no longer the dual of $P^{g,\varepsilon}(t,\tau)$ as demonstrated by \eqref{a6}.
To surmount this impediment, we leverage the regularity properties inherent in  $\mathcal P_4(L^2(\mathcal O))$ and invoke the Vitali theorem to establish the continuity of $\mathcal P^{g,\varepsilon}_*(t,\tau)$ within the restricted domain$(B_{\mathcal P_4(L^2(\mathcal O))},d_{\mathcal P(L^2(\mathcal O))})\times \mathcal H(g_0) $ rather than endeavoring to prove it across the entire space $(\mathcal P_4(L^2(\mathcal O)),d_{\mathcal P(L^2(\mathcal O))})\times \mathcal H(g_0) $.

The first goal of this paper is to prove the existence and uniqueness of uniform measure attractors for the almost periodic external term in  McKean-Vlasov stochastic equation \eqref{a1}-\eqref{a5}
which is dependent on the laws of
solutions.
 To this end, the estimates of the solutions must
be uniform with respect to all translations of the external term involved in the system.

The secondary objective of this paper is to investigate the limiting behavior of the uniform measure attractors associated with the system of \eqref{a1}-\eqref{a5}
 as $\varepsilon\rightarrow0$. Specifically, we aim to understand how these attractors, originating from $(n+1)$-dimensional unbounded thin domains, undergo a collapse into the space $\mathbb R^n$.

  The following sections of the paper are organised as follows: In Section $2$,  we recall the theory of uniform measure attractors for processes defined on the space of probability measures.
  In Section $3$, we reformulate the problem and the transformation from the varying thin domain to the fixed domain, denoted by $\mathcal O$.
   Section $4$ is devoted to the uniform  estimates and the tail estimates of the solutions.
In Section $5$, we show the existence and uniqueness of almost periodic measures of \eqref{a1}-\eqref{a5}.
In the last section, we prove the upper semi-continuity of uniform measure
attractors of \eqref{a1}-\eqref{a5} as $\varepsilon\rightarrow0$.

\section{Uniform measure attractors}

This section reviews the theory of uniform measure attractors for processes in the space of probability measures. Since this space is metrizable, the processes can be seen as in a metric space.

In what follows, we will denote the separable Banach space with norm $\|\cdot\|_X$ by $X$.
Let us define the space of bounded continuous functions on $X$ as $C_b(X)$  equipped with the norm
$$\|\varphi\|_{\infty}=\sup_{x\in X}\vert\varphi(x)\vert.$$
Let $L_b(X)$ denote the space of bounded Lipschitz functions on $X$ which consists of all
functions $\varphi\in C_b(X)$ such that
$$\mathrm{Lip}(\varphi):=\sup_{x_1,x_2\in X, x_1\neq x_2}\frac{|\varphi(x_1)-\varphi(x_2)|}{\|x_1-x_2\|_X}<\infty.$$
The space $L_b(X)$ is equipped with the norm
$$\|\varphi\|_{L_b}=\|\varphi\|_{\infty}+\mathrm{Lip}(\varphi).$$
Denote by $\mathcal P(X)$ be the set of probability measures on $\left(X,\mathcal B(X)\right)$, where
$\mathcal{B}(X)$ is the Borel $\sigma$-algebra of $X$. Given $\varphi\in C_b(X)$ and $\mu\in\mathcal P(X)$,
we write
$$(\varphi,\mu)=\int_X \varphi(x)\mu(dx).$$
Define a metric on $\mathcal P(X)$ by
$$d_{\mathcal P(X)}(\mu_1,\mu_2)=\sup_{\substack{\varphi\in L_b(X)\\ \|\varphi\|_L\le 1}}|(\varphi,\mu_1)-(\varphi,\mu_2)|,\quad \forall \mu_1,\mu_2\in\mathcal P(X).$$
Then $(\mathcal P(X),d_{\mathcal P(X)})$ is a polish space.
Moreover, a sequence $\{\mu_n\}^{\infty}_{n=1}\subset\mathcal P(X)$ convergence to $\mu$ in
$(\mathcal P(X),d_{\mathcal P(X)})$ if and only if $\{\mu_n\}^{\infty}_{n=1}$ convergence to $\mu$ weakly.

Given $p\ge1$, denote by $\left(\mathcal{P}_p(X),\mathbb W_p\right)$ as defined by
$$
\mathcal{P}_p(X)=\left\{\mu \in \mathcal{P}(X): \int_X\|x\|_X^p \mu(d x)<\infty\right\},
$$
and
$$\mathbb W_p(\mu,\nu)=\inf_{\pi\in \Pi(\mu,\nu)}\left(\int_{X\times X}\|x-y\|^p_X(dx,dy)\right)^{\frac{1}{p}},$$
where $\Pi(\mu,\nu)$ is the set of all coupling of $\mu$ and $\nu$. The metric $\mathbb W_p$ is called the Wasserstein distance.

  Given $r>0$, denote by
$$
B_{\mathcal{P}_p(X)}(r)=\left\{\mu \in \mathcal{P}_p(X):\left(\int_X\|x\|_X^p \mu(d x)\right)^{\frac{1}{p}} \leq r\right\}.
$$
A subset $S\subset\mathcal P_p(X)$ is bounded if there is $r>0$ such that $S\subset B_{\mathcal P_p(X)}(r)$. If $S$ is bounded in
$\mathcal P_p(X)$, then we set
$$\|S\|_{\mathcal P_p(X)}=\sup_{\mu\in S}\left(\int_X\|x\|^p_X\mu(dx)\right)^{\frac{1}{p}}.$$
Note that $(\mathcal{P}_p(X),\mathbb W_p)$ is a polish space, but $(\mathcal{P}_p(X),d_{\mathcal P(X)})$ is not complete. Since for every
$r>0$, $B_{\mathcal{P}_p(X)}(r)$ is a closed subset of $\mathcal P(X)$ with respect to the metric $d_{\mathcal P(X)}$, we know that the space
$(B_{\mathcal{P}_p(X)}(r),d_{\mathcal P(X)})$ is complete for every $r>0$.

Recall that the Hausdorff semi-metric between subsets of $\mathcal{P}_p(X)$ is given by

$$
d_{\mathcal{P}(X)}(Y, Z)=\sup _{y \in Y} \inf _{z \in Z} d_{\mathcal{P}(X)}(y, z), \quad Y, Z \subseteq \mathcal{P}_p(X), \quad Y, Z \neq \emptyset .
$$
We assume that $g_0(t)$ is an almost periodic function in $t \in \mathbb{R}$ with values in $X$. Denote by $C_b(\mathbb{R}, X)$ the space of bounded continuous functions on $\mathbb{R}$ with the norm $\|g\|_{C_b(\mathbb{R}, X)}=\sup\limits_{t \in \mathbb{R}}\|g(t)\|_X$ for $g \in C_b(\mathbb{R}, X)$. Since an almost periodic function is bounded and uniformly continuous on $\mathbb{R}$ (see, e.g., \cite{LZ}), it follows that ${g}_0 \in C_b(\mathbb{R}, X)$. Further, by Bochners criterion in \cite{LZ}, whenever ${g}_0: \mathbb{R} \rightarrow X$ is almost periodic, the set of all translations $\left\{ g_0(\cdot+h): h \in \mathbb{R}\right\}$ is precompact in $C_b(\mathbb{R}, X)$. Let $\mathcal{H}\left({g}_0\right)$ be the closure of this set in $C_b(\mathbb{R}, X)$. Then, for any $g \in \mathcal{H}\left({g}_0\right), g$ is almost periodic and $\mathcal{H}(g)=\mathcal{H}\left({g}_0\right)$. For each $h \in \mathbb{R}$, denote by $T(h)$ the translation on $\mathcal{H}\left({g}_0\right)$ with
$T(h)g=g(\cdot+h)$ for all $g\in\mathcal H({g}_0)$. It is evident that $\{T(h)\}_{h\in\mathbb R}$ is a continuous translation group on
$\mathcal H({g}_0)$ that leaves $\mathcal H({g}_0)$ invariant:
$$T(h)\mathcal H({g}_0)=\mathcal H({g}_0),\quad \text{for all }h\in\mathbb R.$$
\begin{df}\label{b2}
A family $U^g=\left\{U^g(t, \tau): t \geq \tau, \tau \in \mathbb{R}\right\}$ of mappings from $\mathcal{P}_p(X)$ to $\mathcal{P}_p(X)$ is called a process on $\mathcal{P}_p(X)$ with time symbol $g \in \mathcal{H}\left({g}_0\right)$, if for all $\tau \in \mathbb{R}$ and $t \geq s \geq \tau$, the following conditions are satisfied:

(a) $U^g(\tau, \tau)=I_{\mathcal{P}_p(X)}$, where $I_{\mathcal{P}_p(X)}$ is the identity operator on $\mathcal{P}_p(X)$;

(b) $U^g(t, \tau)=U^g(t, s) \circ U^g(s, \tau).\\$
The family of process $\left\{U^g(t, \tau)\right\}_{g \in \mathcal{H}\left(g_0\right)}$ are called jointly continuous if it is continuous in both $\mathcal{P}_p(X)$ and $\mathcal{H}\left(g_0\right)$.
\end{df}
It is assumed that the following translation identity holds for the processes $\left\{U^g(t, \tau)\right\}_{g \in \mathcal{H}\left(g_0\right)}$ and the translation group $\{T(h)\}_{h \in \mathbb{R}}$ :
\begin{equation}
U^g(t+h, \tau+h)=U^{T(h) g}(t, \tau), \quad \text { for all } h \in \mathbb{R},\ t \geq \tau \text { and } \tau \in \mathbb{R}.
\end{equation}

\begin{df}
 A closed set $B \subset \mathcal{P}_p(X)$ is called a uniform absorbing set of the family of processes $\left\{U^g(t, \tau)\right\}_{g \in \mathcal{H}\left(g_0\right)}$ with respect to $g \in \mathcal{H}\left(g_0\right)$ if for any bounded $D \subset \mathcal{P}_p(X)$, there exists $T=T\left(D, g_0\right)>0$ such that
$$
U^g(t,0) D \subseteq B, \quad \text { for all } g \in \mathcal{H}\left(g_0\right) \text { and } t \geq T.
$$
\end{df}

\begin{df} The family of processes $\left\{U^g(t, \tau)\right\}_{g \in \mathcal{H}\left(g_0\right)}$ is said to be uniformly asymptotically compact in $\mathcal{P}_p(X)$ with respect to $g \in \mathcal{H}\left(g_0\right)$ if $\left\{U^{g_n}\left(t_n, 0\right) \mu_n\right\}_{n=1}^{\infty}$ has a convergent subsequence in $\mathcal{P}_p(X)$ whenever $t_n \rightarrow+\infty$ and $\left(\mu_n, g_n\right)$ is bounded in $\mathcal{P}_p(X) \times \mathcal{H}\left(g_0\right)$.
\end{df}

\begin{df}
A set $\mathcal{A}$ of $\mathcal{P}_p(X)$ is called a uniform measure attractor of the family of processes $\left\{U^g(t, \tau)\right\}_{g \in \mathcal{H}\left(g_0\right)}$ with respect to $g \in \mathcal{H}\left(g_0\right)$ if the following conditions are satisfied,

(i) $\mathcal{A}$ is compact in $\mathcal{P}_p(X)$;

(ii)$\mathcal A$ is uniformly quasi-invariant, that is, for every $\tau\in\mathbb R$ and $t\ge\tau$,
$$\mathcal A\subseteq \bigcup_{g\in\mathcal H(g_0)}U^g(t,\tau)\mathcal A;$$
(iii) $\mathcal{A}$ attracts every bounded set in $\mathcal{P}_p(X)$ uniformly with respect to $g \in \mathcal{H}\left(g_0\right)$, that is, for any bounded $D \subset \mathcal{P}_p(X)$
$$
\lim _{t \rightarrow \infty} \sup _{g \in \mathcal{H}\left(g_0\right)} d\left(U^g(t, \tau) D, \mathcal{A}\right)=0, \quad \text { for all } \tau \in \mathbb{R};
$$
(iv) $\mathcal{A}$ is minimal among all compact subsets of $\mathcal{P}_p(X)$ satisfying property (iii); that is, if $\mathcal{C}$ is any compact subset of $\mathcal{P}_p(X)$ satisfying property (iii), then $\mathcal{A} \subseteq \mathcal{C}$.
\end{df}

\begin{df}
Given $g \in \mathcal{H}\left(g_0\right)$, a mapping $\chi: \mathbb{R} \rightarrow \mathcal{P}_p(X)$ is called a complete solution of $U^g(t, \tau)$ if for every $t \in \mathbb{R}^{+}$and $\tau \in \mathbb{R}$, the following holds:
$$
U^g(t, \tau)\chi(\tau)=\chi(t).
$$
\end{df}
The kernel of the process $U^g(t, \tau)$ is the collection $\mathcal{K}_g$ of all its bounded complete solutions. The kernel section of the process $U^g(t, \tau)$ at time $s \in \mathbb{R}$ is the set
$$
\mathcal{K}_g(s)=\left\{\xi(s): \xi(\cdot) \in \mathcal{K}_g\right\}.
$$
If the family of processes $\left\{U^g(t, \tau)\right\}_{g \in \mathcal{H}\left(g_0\right)}$ has a uniform measure attractor, then it must be unique. To prove the existence of such a uniform measure attractor, it is convenient to transfer the family of processes to a semigroup of nonlinear operators, and then use the semigroup theory to investigate the uniform measure attractor of the processes. As in \cite{CV}, we define a nonlinear semigroup $\{S(t)\}_{t \geq 0}$ acting on the extended phase space $\mathcal{P}_p(X) \times \mathcal{H}\left(g_0\right)$ by the following formula, for every $t \geq 0, \mu \in \mathcal{P}_p(X)$ and $g \in \mathcal{H}\left(g_0\right)$,
$$
S(t)(\mu, g)=\left(U^g(t, 0) \mu, T(t) g\right).
$$
By the translation identity and Definition \ref{b2} of the process, it is clear that $\{S(t)\}_{t \geq 0}$ satisfies the semigroup identities: for any $t \geq s \geq 0$,
$$
S(0)=\mathrm{I}_{\mathcal{P}_{\mathrm{p}}(\mathrm{X}) \times \mathcal{H}\left({g}_0\right)}, \quad S(t) S(s)=S(t+s).
$$

We know from \cite{CV} that if $\{S(t)\}_{t \geq 0}$ has a global attractor in the extended phase space $\mathcal{P}_p(X) \times \mathcal{H}\left(g_0)\right.$ then the family of the processes $\left\{U^g(t, \tau)\right\}_{g \in \mathcal{H}\left(g_0\right)}$ possesses a uniform measure attractor in the phase space $\mathcal{P}_p(X)$, which is actually the projection onto $\mathcal{P}_p(X)$ of the global attractor of $\{S(t)\}_{t \geq 0}$.

In consequence of the uniform attractors theory set forth in  \cite{CV}, we have the following theorem for the family of processes $\left\{U^g(t, \tau)\right\}_{g \in \mathcal{H}\left(g_0\right)}$. We also refer the reader to \cite{CM02,H88,T12,SY} for the attractors theory of semigroups.

\begin{thm}
 If the semigroup $S(t)$ is continuous, point dissipative and asymptotically compact, then it has a global attractor $\mathcal{A}_S$ in $\mathcal{P}_p(X) \times \mathcal{H}\left(g_0\right)$. Further, if $\mathcal{A}$ is the projection of $\mathcal{A}_S$ onto $\mathcal{P}_p(X)$, then $\mathcal{A}$ is the uniform measure attractor for the family of processes $\left\{U^g(t, \tau)\right\}_{g \in \mathcal{H}\left(g_0\right)}$. In addition,
$$
\mathcal{A}=\underset{g \in \mathcal{H}\left(g_0\right)}{\cup} \mathcal{K}_g(0).
$$
\end{thm}
In accordance with the aforementioned notation from reference \cite{LLZ}, the following criterion is established for the existence and uniqueness of uniform measure attractors.
\begin{thm}\label{b1}
If the family of processes $\left\{U^g(t, \tau)\right\}_{g \in \mathcal{H}\left(g_0\right)}$ is jointly continuous and uniformly asymptotically compact and has a uniform absorbing set $B$, then it has a uniform measure attractor $\mathcal{A}$. In addition,
$$
\mathcal{A}=\underset{g \in \mathcal{H}\left(g_0\right)}{\cup} \mathcal{K}_g(0).
$$
\end{thm}



\section{Existence and uniqueness of solutions}

\setcounter{equation}{0}
In this section, we consider the following equation
\begin{equation}\label{c1}
\left\{\begin{array}{l}
d \hat{u}^{\varepsilon}(t)-\Delta\hat{u}^{\varepsilon}(t)dt+\lambda \hat{u}^{\varepsilon}(t)dt+f\left( x, \hat{u}^{\varepsilon}(t),\mathcal L_{\hat{u}^{\varepsilon}(t)}\right) d t=g(t, x^*) d t \\
\quad+\sum\limits_{k=1}^{\infty}\left(\sigma_{k}(x)+\kappa(x^*)\varpi_k\left(\hat{u}^{\varepsilon}(t),\mathcal L_{\hat{u}^{\varepsilon}(t)}\right)\right) d W_k(t), \quad x\in\mathcal O_{\varepsilon},\ t>\tau \\
\frac{\partial \hat{u}^{\varepsilon}}{\partial \nu_{\varepsilon}}=0,\quad x\in\partial\mathcal O_{\varepsilon},
\end{array}\right.
\end{equation}
with initial data
 \begin{equation}\label{c2}
\hat{u}^{\varepsilon}(\tau,x)=\hat{\xi}^{\varepsilon}(x),\quad x\in\mathcal O_{\varepsilon}.
 \end{equation}
 Throughout this paper, we use $\delta_0$ for the Dirac probability measure at $0$, and we assume $f:\widetilde{\mathcal{O}}\times\mathbb R\times\mathcal P_2(L^2(\mathcal O))\rightarrow\mathbb R$ is continuous and differentiable with respect to the first and second arguments, which further satisfies the conditions:

$\mathbf{(A1)}$. for all $x\in\widetilde{\mathcal{O}}$, $u, u_1,u_2\in\mathbb R$ and $\mu,\mu_1,\mu_2\in\mathcal P_2(L^2(\mathcal O))$,
\begin{align}\label{c3}
f(x,0,\delta_0)=0,
\end{align}
\begin{align}\label{c11}
f(x,u,\mu)u\ge\alpha_1|u|^p-\phi_1(x^*)(1+|u|^2)-\psi_1(x^*)\mu(\|\cdot\|^2),
\end{align}
  \begin{align}\label{c4}
|f(x,u_1,\mu_1)-f(x,u_2,\mu_2)|&\le\alpha_2(\phi_2(x^*)+|u_1|^{p-2}+|u_2|^{p-2})|u_1-u_2|\nonumber\\
&+\phi_3(x^*)\mathbb W_2(\mu_1,\mu_2),
\end{align}
 \begin{align}\label{c14}
\frac{\partial f}{\partial u}(x,u,\mu)\ge-\phi_4(x^*),
\end{align}
\begin{align}\label{c15}
|\frac{\partial f}{\partial x}(x,u,\mu)|\le\phi_5(x^*)\left(1+|u|+\sqrt{\mu(\|\cdot\|^2)}\right),
\end{align}
where $p\ge2$, $\alpha_1,\alpha_2>0$, $\psi_1\in L^1(\mathbb R^n)\cap L^{\infty}(\mathbb R^n)$ and
$\phi_i\in L^{\infty}(\mathbb R)\cap L^1(\mathbb R^n)$ for $i=1,2,3,4,5$.

It follows from \eqref{c3} and \eqref{c4} that for all $u\in\mathbb R, x\in\widetilde{\mathcal{O}}$ and $\mu\in\mathcal P_2(L^2(\mathcal O))$,
\begin{align}
|f(x,u,\mu)|\le\alpha_3|u|^{p-1}+\phi_6(x^*)(1+\sqrt{\mu(\|\cdot\|^2)}),
\end{align}
for some $\alpha_3>0$ and $\phi_6\in L^{\infty}(\mathbb R)\cap L^1(\mathbb R^n)$.

Here, for the diffusion term $\{\sigma_k\}_{k=1}^{\infty}$ and $\{\varpi_k\}^{\infty}_{k=1}$, we assume the following conditions.

$\mathbf{(A2)}$ The function $\sigma=\{\sigma_k\}_{k=1}^{\infty}$ satisfies for all $x=(x^*,x_{n+1})\in\widetilde{\mathcal{O}}$
\begin{equation}\label{c5}
\vert\sigma_k(x)\vert\le\sigma_{1,k}(x^*),
\end{equation}
and
\begin{equation}
\vert\nabla\sigma_k(x)\vert<\sigma_{2,k}(x^*),
\end{equation}
where $\sigma_1=\{\sigma_{1,k}\}^{\infty}_{k=1},\sigma_2=\{\sigma_{2,k}\}^{\infty}_{k=1}\in L^2(\mathbb R^n,l^2)$ with $\\ \sum\limits^{\infty}_{k=1}\left(\|\sigma_{1,k}\|^2_{L^2(\mathbb R^n)}+\|\sigma_{2,k}\|^2_{L^2(\mathbb R^n)}\right)<+\infty$.

$\mathbf{(A3)}$
For each $k\in\mathbb N$, $\varpi_k: \mathbb R\times \mathcal P_2(L^2(\mathcal O))\rightarrow \mathbb R$ is continuous such that for all $u\in\mathbb R$ and
$\mu\in\mathcal P_2(L^2(\mathcal O))$,
\begin{align}\label{c6}
|\varpi_k(u,\mu)|\le\beta_k\left(1+\sqrt{\mu(\|\cdot\|^2)}\right)+\gamma_k|u|,
\end{align}
where $\beta=\{\beta_k\}^{\infty}_{k=1}$ and $\gamma=\{\gamma_k\}^{\infty}_{k=1}$ are nonnegative sequences with $\sum^{\infty}_{k=1}(\beta^2_k+\gamma^2_k)<\infty$.
Furthermore, we assume $\sigma_k(u,\mu)$ is differentiable in $u$ and Lipschitz continuous in both $u$ and $\mu$ uniformly for $t\in\mathbb R$ in the sense
that for all $u_1,u_2\in\mathbb R$ and $\mu_1,\mu_2\in\mathcal P_2(L^2(\mathcal O))$,
\begin{align}\label{c7}
|\varpi_k(u_1,\mu_1)-\varpi_k(u_2,\mu_2)|\le L_{\varpi,k}\left(|u_1-u_2|+\mathbb W_2(\mu_1,\mu_2)\right),
\end{align}
where $L_{\varpi}=\{L_{\varpi,k}\}^{\infty}_{k=1}$ is a sequence of nonnegative numbers such that $\sum^{\infty}_{k=1}L^2_{\varpi,k}<\infty$.

It follows from \eqref{c7} that for all $u\in\mathbb R$ and $\mu\in\mathcal P_2(L^2(\mathcal O))$,
\begin{align}\label{c16}
|\frac{\partial\varpi_k}{\partial u}(u,\mu)|\le L_{\varpi,k}.
\end{align}

Now, we transfer problem \eqref{c1}-\eqref{c2} into the fixed domain $\mathcal O$. To that end,
define $T_{\varepsilon}(x^*,x_{n+1})=(x^*,\frac{x_{n+1}}{\varepsilon \rho(x^*)})$ for $x=(x^*,x_{n+1})\in\mathcal O_{\varepsilon}$.
Let $y=(y^*,y_{n+1})=T_{\varepsilon}(x^*,x_{n+1})$. Then we have
$$
x^*=y^*, \quad x_{n+1}=\varepsilon \rho\left(y^*\right) y_{n+1},
$$
and
$$
\Delta_x \hat{u}(x)=\frac{1}{\rho} \operatorname{div}_y\left(P_{\varepsilon} u(y)\right),
$$
where we denote by $u(y)=\hat{u}(x)$, $\Delta_x$ is the Laplace operator in $x \in \mathcal{O}_{\varepsilon}$, $\operatorname{div}_y$ is the divergence operator in $y \in \mathcal{O}$, and $P_{\varepsilon}$ is the operator given by
$$
P_{\varepsilon} u(y)=\left(\begin{array}{c}
\rho u_{y_1}-\rho_{y_1} y_{n+1} u_{y_{n+1}} \\
\vdots \\
\rho u_{y_n}-\rho_{y_n} y_{n+1} u_{y_{n+1}} \\
-\sum\limits_{i=1}^n y_{n+1} \rho_{y_i} u_{y_i}+\frac{1}{\varepsilon^2 \rho}\left(1+\sum\limits_{i=1}^n\left(\varepsilon y_{n+1} \rho_{y_i}\right)^2\right) u_{y_{n+1}}
\end{array}\right) .
$$

We often write $f( x, u,\mu)$ and $\sigma_k(x)$ as $f\left( x^*, x_{n+1}, u,\mu\right)$ and $\sigma_k\left(x^*, x_{n+1}\right)$ for $x=\left(x^*, x_{n+1}\right)$, respectively. For $y=\left(y^*, y_{n+1}\right) \in \mathcal{O}$ and $t, s \in \mathbb{R}$,  let us define $\sigma_{\varepsilon}=\{\sigma_{k,\varepsilon}\}^{\infty}_{k=1}$ and $f_{\varepsilon}$ as follows:
$$
\begin{array}{rll}
f_{\varepsilon}\left( y^*, y_{n+1}, u,\mu\right)=f\left( y^*, \varepsilon  \rho\left(y^*\right) y_{n+1}, u,\mu\right), & & f_0\left(y^*, u,\mu\right)=f\left( y^*, 0, u,\mu\right),
\end{array}
$$
and
$$
\sigma_{k, \varepsilon}\left(y^*, y_{n+1}\right)=\sigma_k\left(y^*, \varepsilon  \rho\left(y^*\right) y_{n+1}\right), \quad \sigma_{k, 0}\left(y^*\right)=\sigma_k\left(y^*, 0\right).
$$

Denote by $l^2$ the space of square summable sequences of real numbers. For every $u^{\varepsilon}\in L^2(\mathcal O)$, $\mu^{\varepsilon}\in\mathcal P_2(L^2(\mathcal O))$, we define a map $\varpi_{\varepsilon}(u^{\varepsilon},\mu^{\varepsilon}):l^2\rightarrow L^2(\mathcal O)$ by
\begin{align}\label{c17}
\varpi_{\varepsilon}(u^{\varepsilon},\mu^{\varepsilon})(\eta)(y)=\sum^{\infty}_{k=1}(\sigma_{k,\varepsilon}(y)+\kappa(y^*)\varpi_k(u^{\varepsilon}(t),\mu^{\varepsilon}))\eta_k,\quad \forall
\eta=\{\eta\}_{k=1}^{\infty}\in l^2,\ y\in\mathcal{O}.
\end{align}
Similarly, For every $u^0\in L^2(\mathbb R^n)$, $\mu^{0}\in\mathcal P_2(L^2(\mathbb R^n))$, we define  $\varpi_{0}(u^0,\mu^{0}):l^2\rightarrow L^2(\mathbb R^n)$ by
\begin{align}
&\varpi_{0}(u^{0},\mu^{0})(\eta)(y^*)\nonumber\\
&=\sum\limits_{k=1}^{\infty}\left(\sigma_{k,0}(y^*)+\kappa(y^*)\varpi_k\left(u^{0}(t),\mathcal L_{u^{0}(t)}\right)\right)\eta_k,\quad
 \forall
\eta=\{\eta\}_{k=1}^{\infty}\in l^2,\ y^*\in\mathbb R^n.
\end{align}
Then problem \eqref{c1}-\eqref{c2} is equivalent to the following systems for $y=(y^*,y_{n+1})\in\mathcal{O}$,
 \begin{equation}\label{c36}
\left\{\begin{array}{l}
d u^{\varepsilon}(t)-\left(\frac{1}{\rho}\mathrm{div_y}P_{\varepsilon}u^{\varepsilon}(t)-\lambda u^{\varepsilon}(t)\right)dt+f_{\varepsilon}\left( y, u^{\varepsilon}(t),\mathcal L_{u^{\varepsilon}(t)}\right) d t=g(t, y^*) d t \\
\quad+\sum\limits_{k=1}^{\infty}\left(\sigma_{k,\varepsilon}(y)+\kappa(y^*)\varpi_k\left(u^{\varepsilon}(t),\mathcal L_{u^{\varepsilon}(t)}\right)\right) d W_k(t), \quad \ t>\tau, \\
P^{\varepsilon}u^{\varepsilon}\cdot\nu=0,\quad y\in\partial\mathcal O,
\end{array}\right.
\end{equation}
with initial condition
\begin{align}\label{c37}
u^{\varepsilon}(\tau,y)=\xi^{\varepsilon}(y)=\hat{\xi}^{\varepsilon}(T^{-1}_{\varepsilon}(y)),
\end{align}
where $\nu$ is the unit outward normal vector to $\partial\mathcal O$.

Then, we define an inner product $(\cdot,\cdot)_{H_{\rho}(\mathcal O)}$ on $L^2(\mathcal O)$
$$
(u, v)_{H_{\rho}(\mathcal{O})}=\int_{\mathcal{O}}  \rho uvdy, \quad\text {for all } u, v \in L^2(\mathcal{O})
$$
and denote by $L^2(\mathcal{O})$ equipped with this inner product.

For a given value of $0 <\varepsilon \leq 1$, define $a_{\varepsilon}(\cdot, \cdot): H^1(\mathcal{O}) \times H^1(\mathcal{O}) \rightarrow \mathbb{R}$ by
$$
a_{\varepsilon}(u, v)=\left(J^* \nabla_y u, J^* \nabla_y v\right)_{H_{\rho}(\mathcal{O})} \quad \text { for } u, v \in H^1(\mathcal{O}),
$$
where
$$
J^* \nabla_y u(y)=\left(u_{y_1}-\frac{\rho_{y_1}}{\rho} y_{n+1} u_{y_{n+1}}, \ldots, u_{y_n}-\frac{\rho_{y_n}}{\rho} y_{n+1} u_{y_{n+1}}, \frac{1}{\varepsilon \rho} u_{y_{n+1}}\right) .
$$

Define $H_{\varepsilon}^1(\mathcal{O})$ to be the space $H^1(\mathcal{O})$ endowed with the norm
\begin{equation}\label{c19}
\|u\|_{H_{\varepsilon}^1(\mathcal{O})}=\left(\|u\|_{H^1(\mathcal{O})}^2+\frac{1}{\varepsilon^2}\left\|u_{y_{n+1}}\right\|_{L^2(\mathcal{O})}^2\right)^{\frac{1}{2}}.
\end{equation}

It can be shown that there exist positive constants $\eta_1$, $\eta_2$ and $\varepsilon_0$  such that for all $0<\varepsilon<\varepsilon_0$ and $u \in H^1(\mathcal{O})$,
\begin{align}\label{c30}
\eta_1\|u\|_{H_{\varepsilon}^1(\mathcal{O})}^2 \leq a_{\varepsilon}(u, u)+\|u\|_{L^2(\mathcal{O})}^2 \leq \eta_2\|u\|_{H_{\varepsilon}^1(\mathcal{O})}^2 .
\end{align}

Let $A_{\varepsilon}$ be an unbounded operator given by
$$A_{\varepsilon}u=-\frac{1}{\rho}\mathrm{div_y}(P_{\varepsilon}u),\quad u\in D(A_{\varepsilon})=\{u\in H^2(\mathcal O), P_{\varepsilon}u\cdot\nu=0\ \text{on}\
\partial\mathcal O\}.$$
Then we have
\begin{align}
a_{\varepsilon}(u, v)=\left(A_{\varepsilon}u,v\right)_{H_{\rho}(\mathcal O)},\quad \forall u\in D(A_{\varepsilon}),\ \forall v\in H^1(\mathcal O).
\end{align}

In terms of $A_{\varepsilon}$, problem \eqref{c36}-\eqref{c37} is equivalent to
\begin{equation}\label{c8}
\left\{\begin{array}{l}
d u^{\varepsilon}(t)+\left(A_{\varepsilon}u^{\varepsilon}(t)+\lambda u^{\varepsilon}(t)\right)dt+f_{\varepsilon}\left( y, u^{\varepsilon}(t),\mathcal L_{u^{\varepsilon}(t)}\right) d t=g(t, y^*) d t \\
\quad+\varpi_{\varepsilon}\left(u^{\varepsilon}(t),\mathcal L_{u^{\varepsilon}(t)}\right) d W, \quad \ t>\tau, \\
u^{\varepsilon}(\tau)=\xi^{\varepsilon}.
\end{array}\right.
\end{equation}

To reformulate system \eqref{a2}-\eqref{a3}, similarly, we give an inner product $(\cdot,\cdot)_{H_{\rho}(\mathbb R^n)}$ on
$L^2(\mathbb R^n)$
$$
(u, v)_{H_{\rho}(\mathbb R^n)}=\int_{\mathbb R^n}  \rho uvdy^*, \quad\text {for all } u, v \in L^2(\mathbb R^n),
$$
and denote $\left(L^2(\mathbb R^n),(\cdot,\cdot)_{H_{\rho}(\mathbb R^n)}\right)$ by $H_{\rho}(\mathbb R^n)$.

Let $A_0$ be the operator on $H_{\rho}(\mathbb R^n)$ with domain $D(A_0)= H^2(\mathbb R^n)$ as given by
$$A_0 u=-\frac{1}{\rho}\sum^n_{i=1}(\rho u_{y_i})_{y_i},\quad u\in D(A_0).$$
Note that
$$a_0(u,v)=(A_0u,v)_{H_{\rho}(\mathbb R^n)},\quad \forall u\in D(A_0),\ \forall v\in H^1(\mathbb R^n).$$
In terms of $A_0$, system \eqref{a2}-\eqref{a3} is equivalent to
\begin{equation}\label{c9}
\left\{\begin{array}{l}
d u^0(t)+\left(A_{0}u^{0}(t)+\lambda u^{0}(t)\right)dt+f_{0}\left( y^*, u^{0}(t),\mathcal L_{u^{0}(t)}\right) d t=g(t, y^*) dt \\
\quad+\varpi_0\left(u^{0}(t),\mathcal L_{u^{0}(t)}\right) d W(t), \quad \ t>\tau, \\
u^{0}(\tau)=\xi^{0}.
\end{array}\right.
\end{equation}

Under conditions $\mathbf(A1)$-$\mathbf(A3)$, due to the argument presented in  \cite{CB2024}, we establish that for any $\xi^{\varepsilon}\in L^2(\Omega,\mathcal F_{\tau}, L^2(\mathcal O))$,
$g\in\mathcal H(g_0)$,
system \eqref{c8} has a unique solution $u^{\varepsilon}(t,\tau,\xi^{\varepsilon})$ defined on $[\tau,\infty)$.
In particular, $u^{\varepsilon}(t,\tau, \xi^{\varepsilon}),t\ge \tau$, is a continuous $L^2(\mathcal O)$-valued $\mathcal F_t$-adapted
stochastic process such that for every $T>0$,
\begin{align}\label{c13}
u^{\varepsilon}\in &L^2(\Omega,\mathcal C([\tau,\tau+T],L^2(\mathcal O)))\cap L^2(\Omega, L^2((\tau,\tau+T;H^1(\mathcal O)))\nonumber\\
&\cap L^p(\Omega, L^p(\tau,\tau+T;L^p(\mathcal O))),
\end{align}
and for all $t\ge\tau$, $\mathbb P$-almost surely,
\begin{align}
&u^{\varepsilon}(t)+\int^t_0A_{\varepsilon}u^{\varepsilon}(s)ds+\int^t_0\lambda u^{\varepsilon}(s)ds+\int^t_0f_{\varepsilon}\left(\cdot, u^{\varepsilon}(s),\mathcal L_{u^{\varepsilon}(s)}\right) ds\nonumber\\
&=\xi^{\varepsilon}+\int^t_0g(s, \cdot) ds +\int^t_0\varpi_{\varepsilon}\left(u^{\varepsilon}(s),\mathcal L_{u^{\varepsilon}(s)}\right) d W(s),
\end{align}
in $\left(H^1(\mathcal {O})\cap L^p(\mathcal {O})\right)^*$.

Analogously, for any $\xi^0\in L^2(\Omega,\mathcal F_{\tau}, L^2(\mathbb R^n))$, \eqref{c9}  possesses a unique solution
$u^0(t,\tau,\xi^0)$ that is a continuous $L^2(\mathbb R^n)$-valued  stochastic process, adapted to the filtration $ \mathcal F_{\tau}$, and satisfies
\begin{align}
u^{0}\in &L^2(\Omega,\mathcal C([\tau,\tau+T],L^2(\mathbb R^n)))\cap L^2(\Omega,L^2((\tau,\tau+T;H^1(\mathbb R^n)))\nonumber\\
&\cap L^p(\Omega,L^p(\tau,\tau+T;L^p(\mathbb R^n))),
\end{align}
for every $T>0$.

Let $L_2(l^2,H_{\rho}(\mathcal O))$ be the space of Hilbert-Schmidt operators from $l^2$ to $H_{\rho}(\mathcal O)$ with norm $\|\cdot\|_{L^2(l^2,H_{\rho}(\mathcal O))}$. Then
by \eqref{c5} and \eqref{c6} we infer that the operator $\varpi(u,\mu)$ belongs to $L_2(l^2,H_{\rho}(\mathcal O))$ with norm:
\begin{align}\label{c12}
&\|\varpi_{\varepsilon}(u,\mu)\|^2_{L_2(l^2,H_{\rho}(\mathcal O))}=\sum^{\infty}_{k=1}\int_{\mathcal{O}}\rho|\sigma_{k,\varepsilon}(y)+\kappa(y^*)\varpi_k(u(y),\mu)|^2dy\nonumber\\
&\le2\rho_2\|\sigma_{1}\|^2_{L^2(\mathbb R^n,l^2)}+8\rho\|\kappa\|^2\|\beta\|^2_{l^2}(1+\mu(\|\cdot\|^2))+4\|\kappa\|^2_{L^{\infty}(\mathbb R^n)}\|\gamma\|^2_{l^2}\|u\|^2_{H_{\rho}(\mathcal O)}.
\end{align}
Moreover, by \eqref{c7} we see that for all $u_1,u_2\in\mathbb R$ and $\mu_1,\mu_2\in\mathcal P_2(L^2(\mathcal O))$,
\begin{align}\label{c18}
&\|\varpi(u_1,\mu_1)-\varpi(u_2,\mu_2)\|^2_{L_2(l^2,H_{\rho}(\mathcal O)))}\nonumber\\
&=\rho\sum^{\infty}_{k=1}\int_{\mathcal{O}}|\kappa(x^*)|^2|\varpi_k(u_1(x),\mu_1)-\varpi_k(u_2,\mu_2)|^2dx\nonumber\\
&\le 2\rho\|L_{\varpi}\|^2_{l^2}\left(\|\kappa\|^2_{L^{\infty}(\mathbb R^n)}\|u_1-u_2\|^2+\|\kappa\|^2_{L^2(\mathbb R^n)}\mathbb W^2_2(\mu_1,\mu_2)\right).
\end{align}

In the sequel, we also assume that the coefficient  $\lambda$ is sufficiently large, such that
\begin{align}\label{c10}
\lambda>12\|\kappa\|^2\|\beta\|^2_{l^2}+6\|\kappa\|^2_{L^{\infty}(\mathbb R^n)}\|\gamma\|^2_{l^2}+\|\phi_1\|_{L^{\infty}(\mathbb R^n)}+\|\psi_1\|_{L^{1}(\mathbb R^n)}.
\end{align}
It can be deduced from \eqref{c10} that there exists a sufficiently small number $\eta\in(0,1)$ such that
\begin{align}\label{c38}
2\lambda-3\eta>24\|\kappa\|^2\|\beta\|^2_{l^2}+12\|\kappa\|^2_{L^{\infty}(\mathbb R^n)}\|\gamma\|^2_{l^2}+2\|\phi_1\|_{L^{\infty}(\mathbb R^n)}+2\|\psi_1\|_{L^{1}(\mathbb R^n)}.
\end{align}

\section{Priori moment estimates of solutions}

\setcounter{equation}{0}
In this section, we present uniform estimates concerning the solution $u^{\varepsilon}(t)$, which are crucial for demonstrating
the existence and uniqueness of uniform measure attractors.
\begin{lm}\label{d1}
Under $\mathbf(A1)$-$\mathbf(A3)$ and \eqref{c10} hold, then for every $R>0$, there exists $T=T(R)>0$, independent of $\varepsilon$, such that
for any $\tau\in\mathbb R$, $t-\tau\ge T$, and $0<\varepsilon<\varepsilon_0$, the solution $u^{\varepsilon}$ of \eqref{c8} satisfies
$$\mathbb E(\|u^{\varepsilon}(t,\tau,\xi^{\varepsilon})\|^2_{ H_{\rho}(\mathcal O)})\le M_1,$$
and
\begin{align}\label{d37}
\int^t_{\tau}e^{\eta(s-\tau)}\mathbb E\left(\|u^{\varepsilon}(s,\tau,\xi^{\varepsilon})\|^2_{ H^1_{\varepsilon}(\mathcal O)}+\|u^{\varepsilon}(s,\tau,\xi^{\varepsilon})\|^p_{L^p(\mathcal O)}\right)ds<M_1,
\end{align}
where $\mathbb E(\|\xi^{\varepsilon}\|^2_{ H_{\rho}(\mathcal O)})\le R$, and  $M_1$ is constant depending on $\eta, g_0$. In particular,
$M_1$ is independent $\xi^{\varepsilon},\tau, g_{}\in\mathcal H(g_0)$ and $\varepsilon$.
\begin{proof}
By \eqref{c8} and Ito's formula, we have for $t\ge\tau$,
\begin{align}\label{d2}
\begin{split}
&e^{\eta t}\|u^{\varepsilon}(t)\|^2_{ H_{\rho}(\mathcal O)}+2\int^t_{\tau}e^{\eta s}a_{\varepsilon}\left(u^{\varepsilon}(s),u^{\varepsilon}(s)\right)ds
+(2\lambda-\eta)\int^t_{\tau}e^{\eta s}\|u^{\varepsilon}(s)\|^2_{ H_{\rho}(\mathcal O)}ds\\
&+2\int^t_{\tau}e^{\eta s}\big(f_{\varepsilon}(\cdot,u^{\varepsilon}(s),\mathcal L_{u^{\varepsilon}(s)}),u^{\varepsilon}(s)\big)_{ H_{\rho}(\mathcal O)}ds\\
&=e^{\eta \tau}\|\xi^{\varepsilon}\|^2_{ H_{\rho}(\mathcal O)}
+2\int^t_{\tau}e^{\eta s}(g_{}(s,\cdot),u^{\varepsilon}(s))_{ H_{\rho}(\mathcal O)}ds\\
&+\int^t_{\tau}e^{\eta s}\|\varpi_{\varepsilon}(u^{\varepsilon}(s),\mathcal L_{u^{\varepsilon}(s)})\|^2_{ L_2(l^2,H_{\rho})}ds\\
&+2\int^t_{\tau}e^{\eta s}\left(\varpi_{\varepsilon}(u^{\varepsilon}(s),\mathcal L_{u^{\varepsilon}(s)}),u^{\varepsilon}(s)\right)_{ H_{\rho}(\mathcal O)}dW(s),
\end{split}
\end{align}
$\mathbb P$-almost surely. For each $m\in N$, define a stopping time $\tau_m$ as follows:
$$  \tau_m=\inf\{t\ge\tau:\|u^{\varepsilon}(t)\|_{ H_{\rho}(\mathcal O)}>m\}.$$
As is customary, we denote $\inf\empty=+\infty$. Utilizing  \eqref{d2} we can derive the following for all  $t\ge\tau$,
\begin{align}\label{d3}
\begin{split}
&\mathbb E\left(e^{\eta(t\wedge\tau_m)}\|u^{\varepsilon}(t\wedge\tau_m)\|^2_{ H_{\rho}(\mathcal O)}+2\int^{t\wedge\tau_m}_{\tau}e^{\eta s}a_{\varepsilon}(u^{\varepsilon}(s),u^{\varepsilon}(s))ds\right)\\
&=\mathbb E\left(e^{\eta\tau}\|\xi^{\varepsilon}\|^2_{ H_{\rho}(\mathcal O)}\right)+(\eta-2\lambda)\mathbb E\left(\int^{t\wedge\tau_m}_{\tau}e^{\eta s}\|u^{\varepsilon}(s)\|^2_{H_{\rho}(\mathcal O)}ds\right)\\
&-2\mathbb E\left(\int^{t\wedge\tau_m}_{\tau}\int_{\mathcal O}e^{\eta s}\rho f_{\varepsilon}(y,u^{\varepsilon}(s,y),\mathcal L_{u^{\varepsilon}(s)})u^{\varepsilon}(s,y)dyds\right)\\
&+2\mathbb E\left(\int^{t\wedge\tau_m}_{\tau}e^{\eta s}(g(s,\cdot),u^{\varepsilon}(s))_{ H_{\rho}(\mathcal O)}ds\right)\\
&+\mathbb E\left(\int^{t\wedge\tau_m}_{\tau}e^{\eta s}\|\varpi_{\varepsilon}(u^{\varepsilon}(s),\mathcal L_{u^{\varepsilon}(s)})\|^2_{ L_2(l^2,H_{\rho}(\mathcal O))}ds\right).
\end{split}
\end{align}
Next, we derive the uniform estimates for the terms on right-hand of \eqref{d3}. With regard to the third term on the right-hand side of  \eqref{d3}
by \eqref{c11}, we obtain
\begin{align}
\begin{split}
&-2\mathbb E\left(\int^{t\wedge\tau_m}_{\tau}\int_{\mathcal O}e^{\eta s}\rho f_{\varepsilon}(y,u^{\varepsilon}(s,y),\mathcal L_{u^{\varepsilon}(s)})u^{\varepsilon}(s,y)dyds\right)\\
&\le-2\alpha_1\mathbb E\left(\int^{t\wedge\tau_m}_{\tau}e^{\eta s}\rho_1\|u^{\varepsilon}(s)\|^p_{L^p(\mathcal O)}ds\right)\\
&+2\mathbb E\left(\int^{t\wedge\tau_m}_{\tau}e^{\eta s}\|\phi_1\|_{L^{\infty}(\mathbb R^n)}\|u^{\varepsilon}(s)\|^2_{H_{\rho}(\mathcal O)}ds\right)\\
&+2\mathbb E\left(\int^{t\wedge\tau_m}_{\tau}e^{\eta s}\left(\|\phi_1\|_{L^{1}(\mathbb R^n)}+\|\psi_1\|_{L^{1}(\mathbb R^n)}\mathbb E\left(\|u^{\varepsilon}(s)\|^2_{H_{\rho}(\mathcal O)}\right)\right)ds\right)\\
&\le-2\alpha_1\mathbb E\left(\int^{t\wedge\tau_m}_{\tau}e^{\eta s}\rho_1\|u^{\varepsilon}(s)\|^p_{L^p(\mathcal O)}ds\right)+2\mathbb E\left(\int^{t}_{\tau}e^{\eta s}\|\phi_1\|_{L^{\infty}(\mathbb R^n)}\|u^{\varepsilon}(s)\|^2_{H_{\rho}(\mathcal O)}ds\right)\\
&+2\mathbb E\left(\int^{t}_{\tau}e^{\eta s}\left(\|\phi_1\|_{L^{1}(\mathbb R^n)}+\|\psi_1\|_{L^{1}(\mathbb R^n)}\mathbb E\left(\|u^{\varepsilon}(s)\|^2_{H_{\rho}(\mathcal O)}\right)\right)ds\right)\\
&\le -2\alpha_1\mathbb E\left(\int^{t\wedge\tau_m}_{\tau}e^{\eta s}\rho_1\|u^{\varepsilon}(s)\|^p_{L^p(\mathcal O)}ds\right)+
2\int^{t}_{\tau}e^{\eta s}\|\phi_1\|_{L^{1}(\mathbb R^n)}ds\\
&+2\int^{t}_{\tau}e^{\eta s}(\|\phi_1\|_{L^{\infty}(\mathbb R^n)}+\|\psi_1\|_{L^{1}(\mathbb R^n)})\mathbb E\left(\|u^{\varepsilon}(s)\|^2_{H_{\rho}(\mathcal O)}\right)ds.
\end{split}
\end{align}
With regard to the fourth term on the right-hand side of \eqref{d3}, we get
\begin{align}
\begin{split}
&2\mathbb E\left(\int^{t\wedge\tau_m}_{\tau}e^{\eta s}(g(s,\cdot),u^{\varepsilon}(s))_{ H_{\rho}(\mathcal O)}ds\right)\\
&\le\int^{t\wedge\tau_m}_{\tau}\eta e^{\eta s}\mathbb E\left(\|u^{\varepsilon}(s)\|^2_{ H_{\rho}(\mathcal O)}\right)ds+\frac{1}{\eta}\int^{t\wedge\tau_m}_{\tau}\rho_2e^{\eta s}\|g(s)\|^2_{L^2(\mathbb R^n)}ds\\
&\le\int^{t}_{\tau}\eta e^{\eta s}\mathbb E\left(\|u^{\varepsilon}(s)\|^2_{ H_{\rho}(\mathcal O)}\right)ds+\frac{1}{\eta}\int^{t}_{\tau}\rho_2e^{\eta s}\|g(s)\|^2_{L^2(\mathbb R^n)}ds.
\end{split}
\end{align}
For the last term of \eqref{d3}, by \eqref{c12} we have
\begin{align}\label{d4}
\begin{split}
&\mathbb E\left(\int^{t\wedge\tau_m}_{\tau}e^{\eta s}\|\varpi_{\varepsilon}(u^{\varepsilon}(s),\mathcal L_{u^{\varepsilon}(s)})\|^2_{ L_2(l^2,H_{\rho}(\mathcal O))}ds\right)\\
&\le2\mathbb E\left(\int^{t\wedge\tau_m}_{\tau}\rho_2e^{\eta s}\|\sigma_{1}\|^2_{ L_2(\mathbb R^n,l^2)})ds\right)\\
&+8\rho\|\kappa\|^2\|\beta\|^2_{l^2}\mathbb E\left(\int^{t\wedge\tau_m}_{\tau}e^{\eta s}(1+\mathbb E(\|u^{\varepsilon}(s)\|^2_{H_{\rho}(\mathcal O)}))ds\right)\\
&+4\|\kappa\|^2_{L^{\infty}(\mathbb R^n)}\|\gamma\|^2_{l^2}\mathbb E\left(\int^{t\wedge\tau_m}_{\tau}e^{\eta s}\|u^{\varepsilon}(s)\|^2_{H_{\rho}(\mathcal O)}ds\right)\\
&\le2\mathbb E\left(\int^{t}_{\tau}\rho_2e^{\eta s}\|\sigma_{1}\|^2_{ L_2(\mathbb R^n,l^2)})ds\right)\\
&+8\rho\|\kappa\|^2\|\beta\|^2_{l^2}\mathbb E\left(\int^{t}_{\tau}e^{\eta s}(1+\mathbb E(\|u^{\varepsilon}(s)\|^2_{H_{\rho}(\mathcal O)}))ds\right)\\
&+4\|\kappa\|^2_{L^{\infty}(\mathbb R^n)}\|\gamma\|^2_{l^2}\mathbb E\left(\int^{t}_{\tau}e^{\eta s}\|u^{\varepsilon}(s)\|^2_{H_{\rho}(\mathcal O)}ds\right)\\
&\le2\frac{1}{\eta}\rho_2e^{\eta t}\|\sigma_1\|^2_{L^2(\mathbb R^n,l^2)}+8\|\kappa\|^2\|\beta\|^2_{l^2}\frac{1}{\eta}e^{\eta t}\\
&+(8\|\kappa\|^2\|\beta\|^2_{l^2}+4\|\kappa\|^2_{L^{\infty}(\mathbb R^n)}\|\gamma\|^2_{l^2})\int^t_{\tau}e^{\eta s}\mathbb E\left(\|u^{\varepsilon}(s)\|^2_{H_{\rho}(\mathcal O)}\right)ds.
\end{split}
\end{align}
From  \eqref{d3}-\eqref{d4}, it can be deduced that for all $t\ge\tau$,
\begin{align}\label{d5}
\begin{split}
&\mathbb E\left(e^{\eta(t\wedge\tau_m)}\|u^{\varepsilon}(t\wedge\tau_m)\|^2_{ H_{\rho}(\mathcal O)}+2\int^{t\wedge\tau_m}_{\tau}e^{\eta s}a_{\varepsilon}(u^{\varepsilon}(s),u^{\varepsilon}(s))ds\right)\\
&+\mathbb E\left((2\lambda-\eta)\int^{t\wedge\tau_m}_{\tau}e^{\eta s}\|u^{\varepsilon}(s)\|^2_{H_{\rho}(\mathcal O)}ds+
2\alpha_1\rho_1\int^{t\wedge\tau_m}_{\tau}e^{\eta s}\|u^{\varepsilon}(s)\|^p_{L^p(\mathcal O)}ds\right)\\
&\le\mathbb E\left(e^{\eta\tau}\|\xi^{\varepsilon}\|^2_{ H_{\rho}(\mathcal O)}\right)\\
&+2\int^{t}_{\tau}e^{\eta s}(\|\phi_1\|_{L^{\infty}(\mathbb R^n)}+\|\psi_1\|_{L^{1}(\mathbb R^n)})\mathbb E\left(\|u^{\varepsilon}(s)\|^2_{H_{\rho}(\mathcal O)}\right)ds\\
&+\frac{1}{\eta}\int^{t}_{\tau}\rho_2e^{\eta s}\|g(s)\|^2_{L^2(\mathbb R^n)}ds+2\frac{1}{\eta}\rho_2e^{\eta t}\|\sigma_1\|^2_{L^2(\mathbb R^n,l^2)}\\
&+8\|\kappa\|^2\|\beta\|^2_{l^2}\frac{1}{\eta}e^{\eta t}+2\frac{1}{\eta}\|\phi_1\|_{L^1(\mathbb R^n)}e^{\eta t}\\
&+(2\eta+8\|\kappa\|^2\|\beta\|^2_{l^2}+4\|\kappa\|^2_{L^{\infty}(\mathbb R^n)}\|\gamma\|^2_{l^2})\int^t_{\tau}e^{\eta s}\mathbb E\left(\|u^{\varepsilon}(s)\|^2_{H_{\rho}(\mathcal O)}\right)ds.
\end{split}
\end{align}
Taking the limit of \eqref{d5} as $m\rightarrow\infty$, by Fatou's lemma we obtain for all $t\ge\tau$,
\begin{align}\label{d18}
\begin{split}
&\mathbb E\left(e^{\eta t}\|u^{\varepsilon}(t)\|^2_{ H_{\rho}(\mathcal O)}+2\int^{t}_{\tau}e^{\eta s}a_{\varepsilon}(u^{\varepsilon}(s),u^{\varepsilon}(s))ds\right)\\
&+\mathbb E\left(\eta\int^{t}_{\tau}e^{\eta s}\|u^{\varepsilon}(s)\|^2_{H_{\rho}(\mathcal O)}ds+
2\alpha_1\rho_1\int^{t}_{\tau}e^{\eta s}\|u^{\varepsilon}(s)\|^p_{L^p(\mathcal O)}ds\right)\\
&\le\mathbb E\left(e^{\eta\tau}\|\xi^{\varepsilon}\|^2_{ H_{\rho}(\mathcal O)}\right)\\
&+2\int^{t}_{\tau}e^{\eta s}(\|\phi_1\|_{L^{\infty}(\mathbb R^n)}+\|\psi_1\|_{L^{1}(\mathbb R^n)})\mathbb E\left(\|u^{\varepsilon}(s)\|^2_{H_{\rho}(\mathcal O)}\right)ds\\
&+\frac{1}{\eta}\int^{t}_{\tau}\rho_2e^{\eta s}\|g(s)\|^2_{L^2(\mathbb R^n)}ds+2\frac{1}{\eta}\rho_2e^{\eta t}\|\sigma_1\|^2_{L^2(\mathbb R^n,l^2)}\\
&+8\|\kappa\|^2\|\beta\|^2_{l^2}\frac{1}{\eta}e^{\eta t}+2\frac{1}{\eta}\|\phi_1\|_{L^1(\mathbb R^n)}e^{\eta t}\\
&+(2\eta-2\lambda+8\|\kappa\|^2\|\beta\|^2_{l^2}+4\|\kappa\|^2_{L^{\infty}(\mathbb R^n)}\|\gamma\|^2_{l^2})\int^t_{\tau}e^{\eta s}\mathbb E\left(\|u^{\varepsilon}(s)\|^2_{H_{\rho}(\mathcal O)}\right)ds.
\end{split}
\end{align}
By \eqref{c38} and \eqref{d18} we get for all $t\ge\tau$,
\begin{align}\label{d19}
\begin{split}
&\mathbb E\left(\|u^{\varepsilon}(t)\|^2_{ H_{\rho}(\mathcal O)}\right)+2\int^{t}_{\tau}e^{\eta (s-t)}\mathbb E\left(a_{\varepsilon}(u^{\varepsilon}(s),u^{\varepsilon}(s))\right)ds\\
&+\mathbb E\left(\eta\int^{t}_{\tau}e^{\eta (s-t)}\|u^{\varepsilon}(s)\|^2_{H_{\rho}(\mathcal O)}ds+
2\alpha_1\rho_1\int^{t}_{\tau}e^{\eta (s-t)}\|u^{\varepsilon}(s)\|^p_{L^p(\mathcal O)}\right)ds\\
&\le e^{-\eta(t-\tau)}\mathbb E\left(\|\xi^{\varepsilon}\|^2_{ H_{\rho}(\mathcal O)}\right)+2\frac{1}{\eta}\|\phi_1\|_{L^1(\mathbb R^n)}+\frac{1}{\eta}\rho_2\|g_0\|^2_{C_b(\mathbb R,L^2(\mathbb R^n)}\\
&+2\frac{1}{\eta}\rho_2\|\sigma_1\|^2_{L^2(\mathbb R^n,l^2)}+8\|\kappa\|^2\|\beta\|^2_{l^2}\frac{1}{\eta}\\
&=e^{-\eta(t-\tau)}\mathbb E\left(\|\xi^{\varepsilon}\|^2_{ H_{\rho}(\mathcal O)}\right)+C,
\end{split}
\end{align}
where $C=2\frac{1}{\eta}\|\phi_1\|_{L^1(\mathbb R^n)}+\frac{1}{\eta}\rho_2\|g_0\|^2_{C_b(\mathbb R,L^2(\mathbb R^n)}
+2\frac{1}{\eta}\rho_2\|\sigma_1\|^2_{L^2(\mathbb R^n,l^2)}+8\|\kappa\|^2\|\beta\|^2_{l^2}\frac{1}{\eta}$.

Since $\mathbb E\left(\|\xi^{\varepsilon}\|^2_{ H_{\rho}(\mathcal O)}\right)\le R$, we have
$$ e^{-\eta(t-\tau)}\mathbb E\left(\|\xi^{\varepsilon}\|^2_{ H_{\rho}(\mathcal O)}\right)\le e^{-\eta(t-\tau)} R\rightarrow0,\ \text{as}\quad t-\tau\rightarrow\infty,$$
and hence there exists $T=T(R)>0$ such that for all $t-\tau>T$,
$$e^{-\eta(t-\tau)}\mathbb E\left(\|\xi^{\varepsilon}\|^2_{ H_{\rho}(\mathcal O)}\right)\le C,$$
Which along with \eqref{d19} concludes the proof.
\end{proof}
\end{lm}

By Lemma \ref{d1}, we have following uniform estimates.
\begin{cor}\label{d7}
Assume that $\mathbf(A1)$-$\mathbf(A3)$ and \eqref{c10} hold, then for every $R>0$, there exists $T=T(R)>1$, independent of $\varepsilon$, such that
for any $\tau\in\mathbb R$, $t-\tau\ge T$, and $0<\varepsilon<\varepsilon_0$,
 the solution $u^{\varepsilon}$ satisfies
$$\int^t_{t-1}\mathbb E\left(a_{\varepsilon}\left(u^{\varepsilon}(s),u^{\varepsilon}(s)\right)\right)ds\le M_2,$$
where $\mathbb E(\|\xi^{\varepsilon}\|^2_{ H_{\rho}(\mathcal O)})\le R$, and $M_2$ is constant depends on $\eta, g_0$. In particular,
$M_2$ is independent $\xi^{\varepsilon},\tau, g_{}\in\mathcal H(g_0)$ and $\varepsilon$.
\end{cor}

\begin{lm}\label{d38}
Assume that $\mathbf(A1)$-$\mathbf(A3)$ and \eqref{c10} hold, then for every $R>0$, there exists $T=T(R)>1$, independent of $\varepsilon$, such that
for any $\tau\in\mathbb R$, $t\in[\tau,\tau+ T]$, and $0<\varepsilon<\varepsilon_0$,
 the solution $u^{\varepsilon}$ satisfies
$$\int^t_{\tau}\mathbb E\left(\|u^{\varepsilon}(s,\tau,\xi^{\varepsilon})\|^2_{H^1_{\varepsilon}(\mathcal O)}\right)ds\le M_3,$$
where $\mathbb E(\|\xi^{\varepsilon}\|^2_{ H_{\rho}(\mathcal O)})\le R$, and $M_3$ is constant depends on $\eta, g_0$. In particular,
$M_3$ is independent $\xi^{\varepsilon},\tau, g_{}\in\mathcal H(g_0)$ and $\varepsilon$.
\begin{proof}
By \eqref{d37} and $t\in[\tau,\tau+T]$, the desired inequality follows.
\end{proof}
\end{lm}
The following lemma is concerned with the uniform estimates of solutions of \eqref{c9} which is similar to Lemma \ref{d38}.
\begin{lm}\label{d39}
Assume that $\mathbf(A1)$-$\mathbf(A3)$ and \eqref{c10} hold, then for every $R>0$, there exists $T=T(R)>1$, such that
for any $\tau\in\mathbb R$, $t\in[\tau,\tau+ T]$,
 the solution $u^{0}$ satisfies
$$\int^t_{\tau}\mathbb E\left(\|u^{0}(s,\tau,\xi^{0})\|^2_{H^1_{\varepsilon}(\mathbb R^n)}\right)ds\le M_4,$$
where $\mathbb E(\|\xi^{0}\|^2_{ H_{\rho}(\mathbb R^n)})\le R$, and $M_4$ is constant depends on $\eta, g_0$. In particular,
$M_4$ is independent $\xi^{0},\tau, g_{}\in\mathcal H(g_0)$ and $\varepsilon$.
\end{lm}

The following inequality from \cite{LWW17} is useful for deriving the uniform estimates of solution in $L^2(\Omega, H^1_{\varepsilon}(\mathcal O))$.
\begin{lm}\label{d8}
If \eqref{c14}-\eqref{c15} hold, then for $u\in D(A_{\varepsilon})$,
$$
\begin{aligned}
&(f_{\varepsilon}(\cdot,u,\mathcal L_{u(t)}),A_{\varepsilon}u)_{H_{\rho}(\mathcal O)}\\
&\le M_5\left(a_{\varepsilon}(u,u)+\|\phi_5\|^2_{L^2(\mathbb R^n)}(1+\|u(t)\|^2_{H_{\rho}(\mathcal O)}+\mathbb E(\|u(t)\|^2_{H_{\rho}(\mathcal O)})\right),
\end{aligned}
$$
where $M_5$ is a positive constant independent of $\varepsilon$.
\end{lm}
Next, we establish the uniform estimates of solutions of \eqref{c8} in $L^2(\Omega, H^1_{\varepsilon}(\mathcal O))$.
\begin{lm}\label{d10}
Suppose $\mathbf(A1)$-$\mathbf(A3)$ and \eqref{c10} hold, then for every $R>0$, there exists $T=T(R)>1$, independent of $\varepsilon$, such that
for any $\tau\in\mathbb R$, $t-\tau\ge T$, and $0<\varepsilon<\varepsilon_0$, the solution $u^{\varepsilon}$ satisfies
\begin{align}
\mathbb E(\|u^{\varepsilon}(t,\tau,\xi^{\varepsilon})\|^2_{H^{1}_{\varepsilon}(\mathcal O)})\le
M_6,
\end{align}
 where $\mathbb E(\|\xi^{\varepsilon}\|^2_{ H_{\rho}(\mathcal O)})\le R$, and $M_6$ is constant depends on $\lambda, g_0$. In particular,
$M_6$ is independent is independent $\xi^{0},\tau, g_{}\in\mathcal H(g_0)$ and $\varepsilon$.
\begin{proof}
From \eqref{c8} and Ito's formula that for $\tau\in\mathbb R$, $t-\tau>T$ and $\varsigma\in(t-1,t)$,
\begin{align}\label{d20}
\begin{split}
&a_{\varepsilon}(u^{\varepsilon}(t,\tau,\xi^{\varepsilon}),u^{\varepsilon}(t,\tau,\xi^{\varepsilon}))\\
&+2\int_{\varsigma}^t\|A_{\varepsilon}u^{\varepsilon}(s,\tau,\xi^{\varepsilon})\|^2_{H_{\rho}(\mathcal O)}ds+2\lambda\int^t_{\varsigma}a_{\varepsilon}(u^{\varepsilon}(t,\tau,\xi^{\varepsilon}),u^{\varepsilon}(t,\tau,\xi^{\varepsilon}))ds\\
&=a_{\varepsilon}(u^{\varepsilon}(\varsigma,\tau,\xi^{\varepsilon}),u^{\varepsilon}(\varsigma,\tau,\xi^{\varepsilon}))-2\int^t_{\varsigma}(f_{\varepsilon}(\cdot,u^{\varepsilon}
(s,\tau,\xi^{\varepsilon}),\mathcal L_{u^{\varepsilon}(s)})),A_{\varepsilon}u^{\varepsilon}(s,\tau,\xi^{\varepsilon}))_{ H_{\rho}(\mathcal O)}ds\\
&+2\int^t_{\varsigma}(g_{}(s,\cdot),A_{\varepsilon}u^{\varepsilon}(s,\tau,\xi^{\varepsilon}))_{ H_{\rho}(\mathcal O)}ds\\
&+\sum^{\infty}_{k=1}\int^t_{\varsigma}a_{\varepsilon}\left(\sigma_{k,\varepsilon}+\kappa\varpi_k(u^{\varepsilon}(s,\tau,\xi^{\varepsilon}),\mathcal L_{u^{\varepsilon}(s)}),\sigma_{k,\varepsilon}
+\kappa\varpi_k(u^{\varepsilon}(s,\tau,\xi^{\varepsilon}),\mathcal L_{u^{\varepsilon}(s)})\right)ds\\
&+2\int^t_{\varsigma}\left(\varpi_{\varepsilon}(u^{\varepsilon}(s,\tau,\xi^{\varepsilon}),\mathcal L_{u^{\varepsilon}(s)}),A_{\varepsilon}u^{\varepsilon}(s,\tau,\xi^{\varepsilon})\right)_{ H_{\rho}(\mathcal O)}dW(s).
\end{split}
\end{align}
By Lemma \ref{d8} we know
\begin{align}
&2\int^t_{\varsigma}(f_{\varepsilon}(\cdot,u^{\varepsilon}(s,\tau,\xi^{\varepsilon}),\mathcal L_{u^{\varepsilon}(s)}),A_{\varepsilon}u^{\varepsilon}(s,\tau,\xi^{\varepsilon}))_{ H_{\rho}(\mathcal O)}ds\nonumber\\
&\le 2M_5\int^t_{\varsigma}a_{\varepsilon}(u^{\varepsilon}(s,\tau,\xi^{\varepsilon}),u^{\varepsilon}(s,\tau,\xi^{\varepsilon}))ds\nonumber\\
&+2M_5\|\phi_5\|^2_{L^2(\mathbb R^n)}\int^t_{\varsigma}1+\|u^{\varepsilon}(s)\|^2_{H_{\rho}(\mathcal O)}+\mathbb E(\|u^{\varepsilon}(s)\|^2_{H_{\rho}(\mathcal O)})ds.
\end{align}
For the third term on the right-hand side of \eqref{d20} we have
\begin{align}
&2\int^t_{\varsigma}(g_{}(s,\cdot),A_{\varepsilon}u^{\varepsilon}(s,\tau,\xi^{\varepsilon}))_{ H_{\rho}(\mathcal O)}ds\nonumber\\
&\le \int^t_{\varsigma}\|A_{\varepsilon}u^{\varepsilon}(s,\tau,\xi^{\varepsilon})\|^2_{ H_{\rho}(\mathcal O)}ds+\rho_2\int^t_{\varsigma}\|g_{}(s)\|^2_{L^2(\mathbb R^n)}ds.
\end{align}
For the fourth term on the right-hand side of \eqref{d20}, by \eqref{c16}
 we obtain
\begin{align}\label{d6}
&\sum^{\infty}_{k=1}\int^t_{\varsigma}a_{\varepsilon}\left(\sigma_{k,\varepsilon}+\kappa\varpi_k(u^{\varepsilon}(s,\tau,\xi^{\varepsilon}),\mathcal L_{u^{\varepsilon}(s)}),\sigma_{k,\varepsilon}
+\kappa\varpi_k(u^{\varepsilon}(s,\tau,\xi^{\varepsilon}),\mathcal L_{u^{\varepsilon}(s)})\right)ds
\nonumber\\
&\le2\sum^{\infty}_{k=1}\int^t_{\varsigma}a_{\varepsilon}\left(\sigma_{k,\varepsilon},\sigma_{k,\varepsilon}\right)
+a_{\varepsilon}\left(\kappa\varpi_k(u^{\varepsilon}(s,\tau,\xi^{\varepsilon}),\mathcal L_{u^{\varepsilon}(s)}),\kappa\varpi_k(u^{\varepsilon}(s,\tau,\xi^{\varepsilon}),\mathcal L_{u^{\varepsilon}(s)})\right)ds\nonumber\\
&\le2\rho_2\|\sigma_{2}\|^2_{L^2(\mathbb R^n,l^2)}(t-\varsigma)+2\|\kappa\|^2_{L^{\infty}(\mathbb R^n)}\sum^{\infty}_{k=1}L_{\varpi,k}^2\int^t_{\varsigma}
a_{\varepsilon}(u^{\varepsilon}(s),u^{\varepsilon}(s))ds.
\end{align}
By \eqref{d20}-\eqref{d6} we obtain
\begin{align}
\begin{split}
&\mathbb E\left(a_{\varepsilon}(u^{\varepsilon}(t,\tau,\xi^{\varepsilon}),u^{\varepsilon}(t,\tau,\xi^{\varepsilon}))\right)\\
&\le\mathbb E\left(a_{\varepsilon}(u^{\varepsilon}(\varsigma),u^\varepsilon(\varsigma))\right)+c_1\int^{t}_\varsigma\mathbb E(a_{\varepsilon}(u^{\varepsilon}(s),u^\varepsilon(s)))ds\\
&+2M_5\|\phi_5\|^2_{L^2(\mathbb R^n)}(t-\varsigma)+4M_5\|\phi_5\|^2_{L^2(\mathbb R^n)}\int^t_{\varsigma}\mathbb E(\|u^{\varepsilon}(s)\|^2_{H_{\rho}(\mathcal O)})ds\\
&+\rho_2\int^t_{\varsigma}\|g_{}(s)\|^2_{L^2(\mathbb R^n)}ds
+2\rho_2\|\sigma_{2}\|^2_{L^2(\mathbb R^n,l^2)}(t-\varsigma).
\end{split}
\end{align}
where $c_1=2M_5+2\|\kappa\|^2_{L^{\infty}(\mathbb R^n)}\sum\limits^{\infty}_{k=1}L_{\varpi,k}^2$.
Integrating the inequality with respect to $\varsigma$ over $(t-1,t)$, we find
\begin{align}\label{d9}
\begin{split}
&\mathbb E(a_{\varepsilon}(u^{\varepsilon}(t),u^{\varepsilon}(t)))\le(c_1+1)\int^t_{t-1}\mathbb E(a_{\varepsilon}(u^{\varepsilon}(s),u^{\varepsilon}(s)))ds\\
&+2M_5\|\phi_5\|^2_{L^2(\mathbb R^n)}+4M_5\|\phi_5\|^2_{L^2(\mathbb R^n)}\int^t_{t-1}\mathbb E(\|u(s)\|^2_{H_{\rho}(\mathcal O)})ds\\
&+\rho_2\int^t_{t-1}\|g_{}(s)\|^2_{L^2(\mathbb R^n)}ds+2\rho_2\|\sigma_{2}\|^2_{L^2(\mathbb R^n,l^2)}\\
&\le(c_1+1)\int^t_{t-1}\mathbb E(a_{\varepsilon}(u^{\varepsilon}(s),u^{\varepsilon}(s)))ds
+2M_5\|\phi_5\|^2_{L^2(\mathbb R^n)}\\
&+4M_5\|\phi_5\|^2_{L^2(\mathbb R^n)}\int^t_{t-1}\mathbb E(\|u(s)\|^2_{H_{\rho}(\mathcal O)})ds\\
&+\rho_2\|g_{0}\|^2_{C_b(\mathbb R, L^2(\mathbb R^n))}+2\rho_2\|\sigma_{2}\|^2_{L^2(\mathbb R^n,l^2)},
\end{split}
\end{align}
which in conjunction with Lemmas \ref{d7} and \ref{d1}, completes the proof.
\end{proof}
\end{lm}

\begin{lm}\label{d13}
Assume that $\mathbf(A1)$-$\mathbf(A3)$ and \eqref{c10} hold,
 then for any $\delta>0$, $\tau\in\mathbb R$ and $R>0$, there exists a positive integer $N=N(\delta)$ and $T=T(\delta,R)>0$, independent of $\varepsilon$, such that
the solution $u^{\varepsilon}$ of system \eqref{c8} satisfies, for all $t-\tau\ge T$, $n\ge N$ and $0<\varepsilon<\varepsilon_0$,
$$\mathbb E\left(\int_{|y^*|\ge\sqrt{2}n}\int^1_0|u^{\varepsilon}(t,\tau,\xi^{\varepsilon})(y^*,y_{n+1})|^2d_{y_{n+1}}dy^*\right)\le\delta.$$
when $\mathbb E(\|\xi^{\varepsilon}\|^2_{ H_{\rho}(\mathcal O)})\le R$.
\begin{proof}
Let $\theta$ be a cut-off smooth function such that for any $s\in\mathbb R^+$,
Let $\theta$ be a  function satisfying $0\le\theta(s)\le 1$ for $s\ge 0$ and
\begin{align}\label{d36}
\theta(s)=0\quad \text{for $0\le s\le1$};  \qquad\text{$\theta(s)=1$\quad for $s\ge2 $}.
\end{align}
Let $M=\sup\limits_{s\in\mathbb R^+}|\theta'(s)|, n$ be a fixed integer and $\theta_n=\theta(\frac{|y^*|^2}{n^2})$. By \eqref{c8} we get
\begin{align}\label{d21}
\begin{split}
&d \theta_n u^{\varepsilon}(t)+\left(\theta_nA_{\varepsilon}u^{\varepsilon}(t)+\lambda\theta_nu^{\varepsilon}\right)dt
+\theta_nf_{\varepsilon}\left( \cdot, u^{\varepsilon}(t),\mathcal L_{u^{\varepsilon}(t)}\right) d t=\theta_ng_{}(t,\cdot) d t \\
&\quad+\theta_n\varpi_{\varepsilon}\left(u^{\varepsilon}(t),\mathcal L_{u^{\varepsilon}(t)}\right) dW(t), \quad \ t>\tau. \\
\end{split}
\end{align}
By \eqref{d21} and Ito's formula, we have
\begin{align}\label{d11}
\begin{split}
&e^{\eta t}\|\theta_n u^{\varepsilon}(t)\|^2_{H_{\rho}(\mathcal O)}+2\int^t_{\tau}e^{\eta s}a_{\varepsilon}(u^{\varepsilon},\theta^2_n u^{\varepsilon}(s))ds+(2\lambda-\eta)\int^t_{\tau}e^{\eta s}\|\theta_n u^{\varepsilon}(s)\|^2ds\\
&+2\int^t_{\tau}\int_{\mathbb R^n}e^{\eta s}\rho\theta^2_nf_{\varepsilon}(y,u^{\varepsilon}(s,y),\mathcal L_{u^{\varepsilon}(s)})u^{\varepsilon}(s,y)dyds\\
&=e^{\eta\tau}\|\theta_n \xi^{\varepsilon}\|^2_{H_{\rho}(\mathcal O)}+2\int^t_{\tau}(g(s,\cdot),\theta_n^2u(s))_{H_{\rho}(\mathcal O)}ds\\
&+\int^t_{\tau}e^{\eta s}\|\theta_n\varpi_{\varepsilon}(u^{\varepsilon}(s),\mathcal L_{u^{\varepsilon}(s)})\|^2_{L_2(l^2,H_{\rho}(\mathcal O))}ds\\
&\quad+2\int^t_{\tau}e^{\eta s}(\theta_n^2u^{\varepsilon}(s),\varpi_{\varepsilon}(u^{\varepsilon}(s),\mathcal L_{u^{\varepsilon}(s)}))_{H_{\rho}(\mathcal O)}dW(s),
\end{split}
\end{align}
$\mathbb P$-almost surely. Given $m\in\mathbb N$, denote by
$$\tau_m=\inf\{t\ge\tau:\|u^{\varepsilon}(t)\|_{H_{\rho}(\mathcal O)}>m\}.$$
By \eqref{d11} we have for all $t\ge\tau$,
\begin{align}\label{d22}
\begin{split}
&\mathbb E\left(e^{\eta (t\wedge\tau_m)}\|\theta_n u^{\varepsilon}(t\wedge\tau_m)\|^2_{H_{\rho}(\mathcal O)}+2\int^{t\wedge\tau_m}_{\tau}e^{\eta s}a_{\varepsilon}(\theta^2_nu^{\varepsilon}(s), u^{\varepsilon}(s))ds\right)\\
&=\mathbb E\left(e^{\eta\tau}\|\theta_n \xi^{\varepsilon}\|^2_{H_{\rho}(\mathcal O)}\right)+(\eta-2\lambda)\mathbb E\left(\int^{t\wedge\tau_m}_{\tau}e^{\eta s}\|\theta_n u^{\varepsilon}(s)\|^2_{H_{\rho}(\mathcal O)}ds\right)\\
&-2\mathbb E\left(\int^{t\wedge\tau_m}_{\tau}\int_{\mathcal O}e^{\eta s}\rho\theta^2_nf_{\varepsilon}(y,u^{\varepsilon}(s,y),\mathcal L_{u^{\varepsilon}(s)})u^{\varepsilon}(s,y)dyds\right)\\
&\quad+2\mathbb E\left(\int^{t\wedge\tau_m}_{\tau}\left(\theta_n g(s,\cdot),\theta_n u^{\varepsilon}(s)\right)_{H_{\rho}(\mathcal O)}ds\right)\\
&\quad+\mathbb E\left(\int^{t\wedge\tau_m}_{\tau}e^{\eta s}\|\theta_n\varpi_{\varepsilon}(u^{\varepsilon}(s),\mathcal L_{u^{\varepsilon}(s)})\|^2_{L_2(l^2,H_{\rho}(\mathcal O))}ds\right).
\end{split}
\end{align}
Note that
\begin{equation}
\begin{aligned}
& -2 \mathbb E\left(\int^{t\wedge\tau_m}_{\tau}e^{\eta s}a_{\varepsilon}(\theta^2_nu^{\varepsilon}(s), u^{\varepsilon}(s))ds\right) \\
& =-2 \mathbb E\left(\int^{t\wedge\tau_m}_{\tau}e^{\eta s}\int_{\mathcal{O}} \rho J^* \nabla_y u^{\varepsilon}(s,y) \cdot J^* \nabla_y\left(\theta^2\left(\frac{\left|y^*\right|^2}{n^2}\right) u^{\varepsilon}(s,y)\right) \mathrm{d} y\right) \\
& =-2 \mathbb E\left(\int^{t\wedge\tau_m}_{\tau}e^{\eta s}\int_{\mathcal{O}} \rho \theta^2\left(\frac{\left|y^*\right|^2}{n^2}\right)\left|J^* \nabla_y u^{\varepsilon}(s,y)\right|^2 \mathrm{~d} yds\right)\\
&-2 \mathbb E\left(\int^{t\wedge\tau_m}_{\tau}e^{\eta s} \int_{\mathcal{O}} \rho u^{\varepsilon}(s,y) \cdot J^* \nabla_y u^{\varepsilon}(s,y) \cdot J^* \nabla_y\left(\theta^2\left(\frac{\left|y^*\right|^2}{n^2}\right)\right) \mathrm{d} yds\right) \\
& \leq 2 \mathbb E\left(\int^{t\wedge\tau_m}_{\tau}e^{\eta s} \int_{\mathcal{O}}|\rho|\left|u^{\varepsilon}(s,y)\right|\left|J^* \nabla_y u^{\varepsilon}(s,y)\right| \cdot\left|J^* \nabla_y\left(\theta^2\left(\frac{\left|y^*\right|^2}{n^2}\right)\right)\right| \mathrm{d} yds\right)\\
& \leq \frac{c_2}{n} \mathbb E\left(\int^{t\wedge\tau_m}_{\tau}e^{\eta s} \int_{\mathcal{O}}|\rho|\left|u^{\varepsilon}(s,y)\right|\left|J^* \nabla_y u^{\varepsilon}\right|dyds\right) \\
& \leq \frac{c_3}{n}\mathbb E\left(\int^{t\wedge\tau_m}_{\tau}e^{\eta s}\left(\left\|u^{\varepsilon}(s)\right\|_{H_{\rho}(\mathcal{O})}^2+\left\|J^* \nabla_y u^{\varepsilon}(s)\right\|_{H_\rho(\mathcal{O})}^2\right) ds\right)\\
& \leq \frac{c_3}{n}\mathbb E\left(\int^{t\wedge\tau_m}_{\tau}e^{\eta s}\left(\left\|u^{\varepsilon}(s)\right\|_{H_{\rho}(\mathcal{O})}^2+a_{\varepsilon}\left(u^{\varepsilon}(s), u^{\varepsilon}(s)\right)\right)ds\right) \\
& \leq \frac{c_4}{n}\mathbb E\left(\int^{t\wedge\tau_m}_{\tau}e^{\eta s}\left\|u^{\varepsilon}(s)\right\|_{H_{\varepsilon}^1(\mathcal{O})}^2ds\right)\\
&\leq \frac{c_4}{n}\int^{t}_{\tau}e^{\eta s}\mathbb E\left(\left\|u^{\varepsilon}(s)\right\|_{H_{\varepsilon}^1(\mathcal{O})}^2\right)ds,
\end{aligned}
\end{equation}
where $c_2,c_3,c_4$ are positive constants independent of $n$.

For the third term on the right-hand side of \eqref{d22}, by \eqref{c11} have
\begin{align}
&-2\mathbb E\left(\int^{t\wedge\tau_m}_{\tau}\int_{\mathcal O}e^{\eta s}\rho\theta^2_nf_{\varepsilon}(y,u^{\varepsilon}(s,y),\mathcal L_{u^{\varepsilon}(s)})u^{\varepsilon}(s,y)dyds\right)\nonumber\\
&\le-2\alpha_1\mathbb E\left(\int^{t\wedge\tau_m}_{\tau}\rho e^{\eta s}\int_{\mathcal O}\theta^2_n|u^{\varepsilon}(s,y)|^pdyds\right)\nonumber\\
&+2\mathbb E\left(\int^{t\wedge\tau_m}_{\tau} e^{\eta s}\|\phi_1\|_{L^{\infty}(\mathbb R^n)}\|\theta_n u^{\varepsilon}(s)\|^2_{H_{\rho}(\mathcal{O})}ds\right)\nonumber\\
&+2\mathbb E\left(\int^{t\wedge\tau_m}_{\tau} \rho e^{\eta s}\int_{\mathcal O}\theta^2_n|\phi_1(y^*)|dyds\right)
\nonumber\\
&+2\mathbb E\left(\int^{t\wedge\tau_m}_{\tau} e^{\eta s} \mathbb E(\|u^{\varepsilon}(s)\|^2_{H_{\rho}(\mathcal O)}) \int_{\mathcal O}\theta_n^2 \psi_1(y^*)dyds\right)\nonumber\\
&\le-2\alpha_1\mathbb E\left(\int^{t\wedge\tau_m}_{\tau}\rho e^{\eta s}\int_{\mathcal O}\theta^2_n|u^{\varepsilon}(s,y)|^pdyds\right)\nonumber\\
&+2\int^{t}_{\tau} e^{\eta s}\|\phi_1\|_{L^{\infty}(\mathbb R^n)}\mathbb E\left(\|\theta_n u^{\varepsilon}(s)\|^2_{H_{\rho}(\mathcal{O})}\right)ds\nonumber\\
&+2\int^{t}_{\tau} \int_{\mathcal O}\rho e^{\eta s}\theta^2_n|\phi_1(y^*)|dyds
+2\int^{t}_{\tau} e^{\eta s} \mathbb E(\|u^{\varepsilon}(s)\|^2_{H_{\rho}(\mathcal O)})ds \int_{\mathcal O}\theta_n^2 \psi_1(y^*)dy.
\end{align}
For the fourth term on the right-hand side of \eqref{d22}, by Young's inequality  we have
\begin{align}
\begin{split}
&2\mathbb E\left(\int^{t\wedge\tau_m}_{\tau}\left(\theta_n g(s,\cdot),\theta_n u^{\varepsilon}(s)\right)_{H_{\rho}(\mathcal O)}ds\right)\\
&\le \mathbb E\left(\int^{t\wedge\tau_m}_{\tau}\eta e^{\eta s}\|\theta_n u^{\varepsilon}(s)\|^2_{H_{\rho}(\mathcal O)}ds\right)
+\frac{1}{\eta}\mathbb E\left(\int^{t\wedge\tau_m}_{\tau}e^{\eta s}\rho\|\theta_ng(s)\|^2_{L^2(\mathbb R^n)}ds\right)\\
&\le \int^{t}_{\tau}\eta e^{\eta s}\mathbb E\left(\|\theta_n u^{\varepsilon}(s)\|^2_{H_{\rho}(\mathcal O)}\right)ds
+\frac{1}{\eta}\int^{t}_{\tau}e^{\eta s}\rho\|\theta_ng(s)\|^2_{L^2(\mathbb R^n)}ds.
\end{split}
\end{align}
For the last term on the right-hand side of \eqref{d22}, according to \eqref{c12} we have
\begin{align}\label{d23}
\begin{split}
&\mathbb E\left(\int^{t\wedge\tau_m}_{\tau}e^{\eta s}\|\theta_n\varpi_{\varepsilon}(u^{\varepsilon}(s),\mathcal L_{u^{\varepsilon}(s)})\|^2_{L_2(l^2,H_{\rho}(\mathcal O))}ds\right)\\
&\le 2\mathbb E\left(\int^{t\wedge\tau_m}_{\tau}e^{\eta s}\rho\|\theta_n\sigma_{1}\|^2_{L^2(\mathbb R^n,l^2)}ds\right)\\
&+8\|\theta_n\kappa\|^2_{L^2(\mathbb R^n)}\|\beta\|^2_{l^2}\mathbb E\left(\int^{t\wedge\tau_m}_{\tau}e^{\eta s}(1+\mathbb E(\| u^{\varepsilon}(s)\|^2_{ H_{\rho}(\mathcal O)}))ds\right)\\
&+4\|\kappa\|^2_{L^\infty(\mathbb R^n)}\|\gamma\|^2_{l^2}\mathbb E\left(\int^{t\wedge\tau_m}_{\tau}e^{\eta s}\|\theta_n u^{\varepsilon}(s)\|^2_{ H_{\rho}(\mathcal O)}ds\right)\\
&\le 2\int^{t}_{\tau}e^{\eta s}\rho\|\theta_n\sigma_{1}\|^2_{L^2(\mathbb R^n,l^2)}ds+8\|\theta_n\kappa\|^2_{L^2(\mathbb R^n)}\|\beta\|^2_{l^2}\frac{1}{\eta}e^{\eta t}\\
&+8\|\theta_n\kappa\|^2_{L^2(\mathbb R^n)}\|\beta\|^2_{l^2}\int^{t}_{\tau}e^{\eta s}\mathbb E(\| u^{\varepsilon}(s)\|^2_{ H_{\rho}(\mathcal O)})ds\\
&+4\|\kappa\|^2_{L^\infty(\mathbb R^n)}\|\gamma\|^2_{l^2}\int^{t}_{\tau}e^{\eta s}\mathbb E\left(\|\theta_n u^{\varepsilon}(s)\|^2_{ H_{\rho}(\mathcal O)}\right)ds.
\end{split}
\end{align}
It follows from \eqref{d22}-\eqref{d23} that for all $t\ge\tau$,
\begin{align}\label{d24}
\begin{split}
&\mathbb E\left(e^{\eta (t\wedge\tau_m)}\|\theta_n u^{\varepsilon}(t\wedge\tau_m)\|^2_{H_{\rho}(\mathcal O)}+(2\lambda-\eta)\int^{t\wedge\tau_m}_{\tau}e^{\eta s}\|\theta_n u^{\varepsilon}(s)\|^2_{H_{\rho}(\mathcal O)}ds\right)\\
&\le \mathbb E\left(e^{\eta\tau}\|\theta_n \xi^{\varepsilon}\|^2_{H_{\rho}(\mathcal O)}\right)+2\int^{t}_{\tau} \rho e^{\eta s}\|\theta^2_n\phi_1\|_{L^1(\mathbb R^n)}dyds\\
&+\frac{c_4}{n}\int^{t\wedge\tau_m}_{\tau}e^{\eta s}\mathbb E\left(\left\|u^{\varepsilon}\right\|_{H_{\varepsilon}^1(\mathcal{O})}^2\right)ds
+2\int^{t}_{\tau}e^{\eta s}\rho\|\theta_n\sigma_{1}\|^2_{L^2(\mathbb R^n,l^2)}ds\\
&+\frac{1}{\eta}\int^{t}_{\tau}e^{\eta s}\rho\|\theta_ng(s)\|^2_{L^2(\mathbb R^n)}ds
+8\|\theta_n\kappa\|^2_{L^2(\mathbb R^n)}\|\beta\|^2_{l^2}\frac{1}{\eta}e^{\eta t}\\
&+\int^{t}_{\tau} e^{\eta s}(\eta+4\|\kappa\|^2_{L^\infty(\mathbb R^n)}\|\gamma\|^2_{l^2}+2\|\phi_1\|_{L^{\infty}(\mathbb R^n)}
)\mathbb E\left(\|\theta_n u^{\varepsilon}(s)\|^2_{H_{\rho}(\mathcal O)}\right)ds\\
&+
\int^{t}_{\tau} e^{\eta s}(8\|\theta_n\kappa\|^2_{L^2(\mathbb R^n)}\|\beta\|^2_{l^2}+2\|\theta_n^2 \psi_1\|_{L^1(\mathbb R^n)}) \mathbb E(\|u^{\varepsilon}(s)\|^2_{H_{\rho}(\mathcal O)})ds.
\end{split}
\end{align}
Taking the limit of \eqref{d24} as $m\rightarrow\infty$, by Fatou's Lemma we obtain for all $t\ge\tau$,
\begin{align}\label{d25}
\begin{split}
&\mathbb E\left(e^{\eta t}\|\theta_n u^{\varepsilon}(t)\|^2_{H_{\rho}(\mathcal O)}+(2\lambda-\eta)\int^{t}_{\tau}e^{\eta s}\|\theta_n u^{\varepsilon}(s)\|^2_{H_{\rho}(\mathcal O)}ds\right)\\
&\le \mathbb E\left(e^{\eta\tau}\|\theta_n \xi^{\varepsilon}\|^2_{H_{\rho}(\mathcal O)}\right)+2\int^{t}_{\tau} \rho e^{\eta s}\|\theta^2_n\phi_1\|_{L^1(\mathbb R^n)}dyds\\
&+\frac{c_4}{n}\int^{t}_{\tau}e^{\eta s}\mathbb E\left(\left\|u^{\varepsilon}\right\|_{H_{\varepsilon}^1(\mathcal{O})}^2\right)ds
+2\int^{t}_{\tau}e^{\eta s}\rho\|\theta_n\sigma_{1}\|^2_{L^2(\mathbb R^n,l^2)}ds\\
&+\frac{1}{\eta}\int^{t}_{\tau}e^{\eta s}\rho\|\theta_ng(s)\|^2_{L^2(\mathbb R^n)}ds
+8\|\theta_n\kappa\|^2_{L^2(\mathbb R^n)}\|\beta\|^2_{l^2}\frac{1}{\eta}e^{\eta t}\\
&+\int^{t}_{\tau} e^{\eta s}(\eta+4\|\kappa\|^2_{L^\infty(\mathbb R^n)}\|\gamma\|^2_{l^2}+2\|\phi_1\|_{L^{\infty}(\mathbb R^n)}
)\mathbb E\left(\|\theta_n u^{\varepsilon}(s)\|^2_{H_{\rho}(\mathcal O)}\right)ds\\
&+
\int^{t}_{\tau} e^{\eta s}(8\|\theta_n\kappa\|^2_{L^2(\mathbb R^n)}\|\beta\|^2_{l^2}+2\|\theta_n^2 \psi_1\|_{L^1(\mathbb R^n)}) \mathbb E(\|u^{\varepsilon}(s)\|^2_{H_{\rho}(\mathcal O)})ds.
\end{split}
\end{align}
By \eqref{d25} and \eqref{c38} we get for all $t\ge\tau$,
\begin{align}\label{d26}
\begin{split}
&\mathbb E\left(\|\theta_n u^{\varepsilon}(t)\|^2_{H_{\rho}(\mathcal O)}\right)\\
&\le e^{-\eta(t-\tau)}\mathbb E\left(\|\theta_n \xi^{\varepsilon}\|^2_{H_{\rho}(\mathcal O)}\right)+2\rho_2\frac{1}{\eta}\|\theta^2_n\phi_1\|_{L^1(\mathbb R^n)}\\
&+\frac{c_4}{n}\int^{t}_{\tau}e^{\eta (s-t)}\mathbb E\left(\left\|u^{\varepsilon}\right\|_{H_{\varepsilon}^1(\mathcal{O})}^2\right)ds
+2\rho_2\frac{1}{\eta}\|\theta_n\sigma_{1}\|^2_{L^2(\mathbb R^n,l^2)}\\
&+\frac{1}{\eta}\rho_2\|\theta_ng_0\|^2_{C_b(\mathbb R,L^2(\mathbb R^n))}
+8\|\theta_n\kappa\|^2_{L^2(\mathbb R^n)}\|\beta\|^2_{l^2}\frac{1}{\eta}\\
&+
(8\|\theta_n\kappa\|^2_{L^2(\mathbb R^n)}\|\beta\|^2_{l^2}+2\|\theta_n^2 \psi_1\|_{L^1(\mathbb R^n)})\int^{t}_{\tau} e^{-\eta (t-s)} \mathbb E(\|u^{\varepsilon}(s)\|^2_{H_{\rho}(\mathcal O)})ds.
\end{split}
\end{align}
By \eqref{d26} and Lemma \ref{d1} we find that there for exists $c_5>0$ and $T_1(R)\ge1$ such that for all $t-\tau\ge T_1$
\begin{align}\label{d27}
\begin{split}
&\mathbb E\left(\|\theta_n u^{\varepsilon}(t)\|^2_{H_{\rho}(\mathcal O)}\right)\\
&\le e^{-\eta(t-\tau)}\mathbb E\left(\| \xi^{\varepsilon}\|^2_{H_{\rho}(\mathcal O)}\right)+2\rho_2\frac{1}{\eta}\|\theta^2_n\phi_1\|_{L^1(\mathbb R^n)}\\
&+\frac{1}{\eta^2}\rho_2\|\theta_ng_0\|^2_{C_b(\mathbb R,L^2(\mathbb R^n)}
+2\rho_2\frac{1}{\eta}\|\theta_n\sigma_{\varepsilon}\|^2_{L^2(\mathcal O,l^2))}\\
&+\frac{c_4c_5}{ n}
+8\|\theta_n\kappa\|^2_{L^2(\mathbb R^n)}\|\beta\|^2_{l^2}\frac{1}{\eta}\\
&+
(8\|\theta_n\kappa\|^2_{L^2(\mathbb R^n)}\|\beta\|^2_{l^2}+2\|\theta_n^2 \psi_1\|_{L^1(\mathbb R^n)})c_5.
\end{split}
\end{align}
Since $\mathbb E(\|\xi^{\varepsilon}\|^2_{ H_{\rho}(\mathcal O)})\le R$, we have
$$\lim_{t\rightarrow\infty}e^{\eta(\tau-t)}\mathbb E\left(\| \xi^{\varepsilon}\|^2_{H_{\rho}(\mathcal O)}\right)\le\lim_{t\rightarrow\infty}e^{\eta(\tau-t)} R=0,$$
and hence for every $\delta>0$, there exists $T_2=T_2(R,\delta)\ge T_1$ such that for all $t-\tau\ge T_2$
$$e^{\eta(\tau-t)}\mathbb E\left(\|\xi^{\varepsilon}\|^2_{H_{\rho}(\mathcal O)}\right)<\frac{1}{3\delta}.$$
By \eqref{c5}, $\phi_1\in L^1(\mathbb R^n)$, and $g_0\in\mathcal H(g_0)$,
 we find that for every $\delta>0$, there exists $N_1=N_1(\delta)\in\mathbb N$
such that for all $n\ge N_1$,
\begin{align}
\begin{split}
&2\rho_2\frac{1}{\eta}\|\theta^2_n\phi_1\|_{L^1(\mathbb R^n)}+2\rho_2\frac{1}{\eta}\|\theta_n\sigma_{1}\|^2_{L^2(\mathbb R^n,l^2))}
+\frac{1}{\eta^2}\rho_2\|\theta_ng_0\|^2_{C_b(\mathbb R,L^2(\mathbb R^n))}\\
&\le2\rho_2\frac{1}{\eta}\|\int_{|y^*|\ge n}|\phi_1(y^*)|dy+2\rho_2\frac{1}{\eta}\int_{|y^*|\ge n}\sum^{\infty}_{k=1}|\sigma_{1,k}(y^*)|^2dy^*\\
&+\frac{1}{\eta^2}\rho_2\int_{|y^*|\ge n}\sup_{t\in\mathbb R}|g_0(\cdot,y^*)|^2dy\\
&\le\frac{1}{9\delta}+\frac{1}{9\delta}+\frac{1}{9\delta}=\frac{1}{3\delta}.
\end{split}
\end{align}
Note that $\kappa\in L^2(\mathbb R^n)$ and $\psi_{1}\in L^1(\mathbb R^n)$. Thus, there exists $N_2=N_2(\delta)\ge N_1$ such that
for all $n\ge N_2$,
\begin{align}\label{d28}
\begin{split}
&\frac{c_4c_5}{ n}+8\|\theta_n\kappa\|^2_{L^2(\mathbb R^n)}\|\beta\|^2_{l^2}\frac{1}{\eta}+ (8\|\theta_n\kappa\|^2_{L^2(\mathbb R^n)}\|\beta\|^2_{l^2}+2\|\theta_n^2 \psi_1\|_{L^1(\mathbb R^n)})c_5\\
&\le\frac{c_4c_5}{ n}+8\|\beta\|^2_{l^2}\frac{1}{\eta}\int_{|y^*|\ge n}|\kappa(y^*)|^2dy^*\\
&+ \left(8\|\beta\|^2_{l^2}\int_{|y^*|\ge n}|\kappa(y^*)|^2dy^*2+\int_{|y^*|\ge n}|\psi(y^*)|dy^*\right)c_5\\
&\le\frac{1}{9\delta}+\frac{1}{9\delta}+\frac{1}{9\delta}=\frac{1}{3\delta}.
\end{split}
\end{align}
It follows from \eqref{d27}-\eqref{d28} that for all $t-\tau\ge T_2$ and $n\ge N_2$,
$$\mathbb E\left(\int_{|y^*|\ge\sqrt{2}n}\int^1_0|u^{\varepsilon}(t,\tau,\xi^{\varepsilon})(y^*,y_{n+1})|^2d_{y_{n+1}}dy^*\right)\le\frac{1}{\rho_1}\mathbb E\left(\|\theta_n u^{\varepsilon}(t)\|^2_{H_{\rho}(\mathcal O)}\right)\le\frac{\delta}{\rho_1}.$$
This completes the proof.
\end{proof}
\end{lm}
Next, we derive the uniform estimates of solutions of \eqref{c8} in $L^4(\Omega,\mathcal F, H_{\rho}(\mathcal O))$.
\begin{lm}\label{d35}
Under $\mathbf(A1)$-$\mathbf(A3)$ and \eqref{c10} hold, then for every $R>0$, there exists $T=T(R)>0$, independent of $\varepsilon$, such that
for any $\tau\in\mathbb R$, $t-\tau\ge T$, and $0<\varepsilon<\varepsilon_0$, the solution $u^{\varepsilon}$ of \eqref{c8} satisfies
$$\mathbb E(\|u^{\varepsilon}(t,\tau,\xi^{\varepsilon})\|^4_{ H_{\rho}(\mathcal O)})\le M_7,$$
where $\mathbb E(\|\xi^{\varepsilon}\|^4_{ H_{\rho}(\mathcal O)})\le R$, and  $M_7$ is constant depending on $\lambda, g_0$. In particular,
$M_7$ is independent $\xi^{\varepsilon},\tau, g_{}\in\mathcal H(g_0)$ and $\varepsilon$.
\begin{proof}
By \eqref{d2} and Ito's formula, we get for all $t\ge\tau$
\begin{align}\label{d40}
\begin{split}
&e^{2\eta t}\|u^{\varepsilon}(t)\|^4_{ H_{\rho}(\mathcal O)}+4\int^t_{\tau}e^{2\eta s}\|u^{\varepsilon}(s)\|^2_{H_{\rho}(\mathcal O)}a_{\varepsilon}\left(u^{\varepsilon}(s),u^{\varepsilon}(s)\right)ds\\
&\quad+2(2\lambda-\eta)\int^t_{\tau}e^{2\eta s}\|u^{\varepsilon}(s)\|^4_{ H_{\rho}(\mathcal O)}ds\\
&+4\int^t_{\tau}e^{2\eta s}\|u^{\varepsilon}(s)\|^2_{ H_{\rho}(\mathcal O)}\int_{\mathcal O}\rho f_{\varepsilon}(y,u^{\varepsilon}(s,y),\mathcal L_{u^{\varepsilon}(s)})u^{\varepsilon}(s,y)dyds\\
&=e^{2\eta \tau}\|\xi^{\varepsilon}\|^4_{ H_{\rho}(\mathcal O)}+2\int^t_{\tau}e^{2\eta s}\|u^{\varepsilon}(s)\|^2_{ H_{\rho}(\mathcal O)}\|\varpi_{\varepsilon}(u^{\varepsilon}(s),\mathcal L_{u^{\varepsilon}(s)})\|^2_{ L_2(l^2,H_{\rho}(\mathcal O))}ds\\
&+4\int^t_{\tau}e^{2\eta s}\|u^{\varepsilon}(s)\|^2_{ H_{\rho}(\mathcal O)}(g_{}(s,\cdot),u^{\varepsilon}(s))_{ H_{\rho}(\mathcal O)}ds\\
&+4\int^t_{\tau}e^{2\eta s}\rho^2\|u^{\varepsilon,*}(s)\varpi_{\varepsilon}(u^{\varepsilon}(s),\mathcal L_{u^{\varepsilon}(s)})\|^2_{ L_2(l^2,\mathbb R)}ds\\
&+4\int^t_{\tau}e^{2\eta s}\|u^{\varepsilon}(s)\|^2_{ H_{\rho}(\mathcal O)}\left(\varpi_{\varepsilon}(u^{\varepsilon}(s),\mathcal L_{u^{\varepsilon}(s)}),u^{\varepsilon}(s)\right)_{ H_{\rho}(\mathcal O)}dW(s),
\end{split}
\end{align}
$\mathbb P$-almost surely, where $u^{\varepsilon,*}(s)$ is the element in $L^{2,*}(\mathcal O)$ identified by Riesz representation theorem.
Let $\tau_m=\inf\{t\ge\tau:\|u^{\varepsilon}(t)\|_{ H_{\rho}(\mathcal O)}>m\}$. By \eqref{d40} we have for all $t\ge\tau$,
\begin{align}\label{d29}
\begin{split}
&e^{2\eta (t\wedge\tau_m)}\|u^{\varepsilon}(t\wedge\tau_m)\|^4_{ H_{\rho}(\mathcal O)}+4\int^{t\wedge\tau_m}_{\tau}e^{2\eta s}\|u^{\varepsilon}(s)\|^2_{H_{\rho}(\mathcal O)}a_{\varepsilon}\left(u^{\varepsilon}(s),u^{\varepsilon}(s)\right)ds\\
&\quad+2(2\lambda-\eta)\int^{t\wedge\tau_m}_{\tau}e^{2\eta s}\|u^{\varepsilon}(s)\|^4_{ H_{\rho}(\mathcal O)}ds\\
&+4\int^{t\wedge\tau_m}_{\tau}e^{2\eta s}\|u^{\varepsilon}(s)\|^2_{ H_{\rho}(\mathcal O)}\int_{\mathcal O}\rho f_{\varepsilon}(y,u^{\varepsilon}(s,y),\mathcal L_{u^{\varepsilon}(s)})u^{\varepsilon}(s,y)dyds\\
&=e^{2\eta \tau}\|\xi^{\varepsilon}\|^4_{ H_{\rho}(\mathcal O)}+2\int^{t\wedge\tau_m}_{\tau}e^{2\eta s}\|u^{\varepsilon}(s)\|^2_{ H_{\rho}(\mathcal O)}\|\varpi_{\varepsilon}(u^{\varepsilon}(s),\mathcal L_{u^{\varepsilon}(s)})\|^2_{ L_2(l^2,H_{\rho}(\mathcal O))}ds\\
&+4\int^{t\wedge\tau_m}_{\tau}e^{2\eta s}\|u^{\varepsilon}(s)\|^2_{ H_{\rho}(\mathcal O)}(g_{}(s,\cdot),u^{\varepsilon}(s))_{ H_{\rho}(\mathcal O)}ds\\
&+4\int^{t\wedge\tau_m}_{\tau}e^{2\eta s}\rho^2\|u^{\varepsilon,*}(s)\varpi_{\varepsilon}(u^{\varepsilon}(s),\mathcal L_{u^{\varepsilon}(s)})\|^2_{ L_2(l^2,\mathbb R)}ds\\
&+4\int^{t\wedge\tau_m}_{\tau}e^{2\eta s}\|u^{\varepsilon}(s)\|^2_{ H_{\rho}(\mathcal O)}\left(\varpi_{\varepsilon}(u^{\varepsilon}(s),\mathcal L_{u^{\varepsilon}(s)}),u^{\varepsilon}(s)\right)_{ H_{\rho}(\mathcal O)}dW(s).
\end{split}
\end{align}
Taking the expectation of \eqref{d29} we get, for all $t\ge\tau$,
\begin{align}\label{d30}
\begin{split}
&\mathbb E\left(e^{2\eta (t\wedge\tau_m)}\|u^{\varepsilon}(t\wedge\tau_m)\|^4_{ H_{\rho}(\mathcal O)}\right)+2(2\lambda-\eta)\mathbb E\left(\int^{t\wedge\tau_m}_{\tau}e^{2\eta s}\|u^{\varepsilon}(s)\|^4_{ H_{\rho}(\mathcal O)}ds\right)\\
&\le e^{2\eta \tau}\|\xi^{\varepsilon}\|^4_{ H_{\rho}(\mathcal O)}\\
&-4\mathbb E\left(\int^{t\wedge\tau_m}_{\tau}e^{2\eta s}\|u^{\varepsilon}(s)\|^2_{ H_{\rho}(\mathcal O)}\int_{\mathcal O}\rho f_{\varepsilon}(y,u^{\varepsilon}(s,y),\mathcal L_{u^{\varepsilon}(s)})u^{\varepsilon}(s,y)dyds\right)\\
&+4\mathbb E\left(\int^{t\wedge\tau_m}_{\tau}e^{2\eta s}\|u^{\varepsilon}(s)\|^2_{ H_{\rho}(\mathcal O)}(g_{}(s,\cdot),u^{\varepsilon}(s))_{ H_{\rho}(\mathcal O)}ds\right)\\
&+6\mathbb E\left(\int^{t\wedge\tau_m}_{\tau}e^{2\eta s}\|u^{\varepsilon}(s)\|^2_{ H_{\rho}(\mathcal O)}\|\varpi_{\varepsilon}(u^{\varepsilon}(s),\mathcal L_{u^{\varepsilon}(s)})\|^2_{ L_2(l^2,H_{\rho}(\mathcal O))}ds\right).
\end{split}
\end{align}
For the second term on the right-hand side of \eqref{d30}, by \eqref{c11}, we have
\begin{align}
\begin{split}
&-4\mathbb E\left(\int^{t\wedge\tau_m}_{\tau}e^{2\eta s}\|u^{\varepsilon}(s)\|^2_{ H_{\rho}(\mathcal O)}\int_{\mathcal O}\rho f_{\varepsilon}(y,u^{\varepsilon}(s,y),\mathcal L_{u^{\varepsilon}(s)})u^{\varepsilon}(s,y)dyds\right)\\
&\le-4\alpha_1\mathbb E\left(\int^{t\wedge\tau_m}_{\tau} \rho e^{2\eta s}\|u^{\varepsilon}(s)\|^2_{H_{\rho}(\mathcal O)}\|u^{\varepsilon}(s)\|^p_{L^p(\mathcal O)}ds\right)\\
&+4\mathbb E\left(\int^{t\wedge\tau_m}_{\tau}e^{2\eta s}\|\phi_1\|_{L^{\infty}(\mathbb R^n)}\|u^{\varepsilon}(s)\|^4_{H_{\rho}(\mathcal O)}ds\right)\\
&+4\mathbb E\left(\int^{t\wedge\tau_m}_{\tau}e^{2\eta s}\|u^{\varepsilon}(s)\|^2_{H_{\rho}(\mathcal O)}(\rho\|\phi_1\|_{L^1(\mathbb R^n)}+\|\psi_1\|_{L^1(\mathbb R^n)}\mathbb E\left(\|u^{\varepsilon}(s)\|^2_{ H_{\rho}(\mathcal O)}\right)ds\right)\\
&\le 4\int^{t}_{\tau}e^{2\eta s}\|\phi_1\|_{L^{\infty}(\mathbb R^n)}\mathbb E\left(\|u^{\varepsilon}(s)\|^4_{H_{\rho}(\mathcal O)}\right)ds\\
&+4\int^{t}_{\tau}e^{2\eta s}\rho\|\phi_1\|_{L^1(\mathbb R^n)}\mathbb E\left(\|u^{\varepsilon}(s)\|^2_{H_{\rho}(\mathcal O)}\right)ds\\
&+
4\int^{t}_{\tau}e^{2\eta s}\|\psi_1\|_{L^1(\mathbb R^n)}\mathbb E\left(\|u^{\varepsilon}(s)\|^2_{ H_{\rho}(\mathcal O)}\right)^2ds\\
&\le\int^{t}_{\tau}e^{2\eta s}\left(4\|\phi_1\|_{L^{\infty}(\mathbb R^n)}+4\|\psi_1\|_{L^1(\mathbb R^n)}+\frac{1}{3}\eta\right)\mathbb E\left(\|u^{\varepsilon}(s)\|^4_{H_{\rho}(\mathcal O)}\right)ds\\
&+12\frac{1}{\eta}\rho^2\int^{t}_{\tau}e^{2\eta s}\|\phi_1\|^2_{L^1(\mathbb R^n)}ds.
\end{split}
\end{align}
For the third term on the right-hand side of \eqref{d30}, we get
\begin{align}
\begin{split}
&4\mathbb E\left(\int^{t\wedge\tau_m}_{\tau}e^{2\eta s}\|u^{\varepsilon}(s)\|^2_{ H_{\rho}(\mathcal O)}(g_{}(s,\cdot),u^{\varepsilon}(s))_{ H_{\rho}(\mathcal O)}ds\right)\\
&\le2\eta\mathbb E\left(\int^{t\wedge\tau_m}_{\tau}e^{2\eta s}\|u^{\varepsilon}(s)\|^4_{ H_{\rho}(\mathcal O)}ds\right)+\frac{27}{8\eta^3}\rho\int^{t\wedge\tau_m}_{\tau}e^{2\eta s}\|g(s)\|^2_{L^2(\mathbb R^n)}ds\\
&\le2\eta\int^{t}_{\tau}e^{2\eta s}\mathbb E\left(\|u^{\varepsilon}(s)\|^4_{ H_{\rho}(\mathcal O)}\right)ds+\frac{27}{8\eta^3}\rho\int^{t}_{\tau}e^{2\eta s}\|g(s)\|^2_{L^2(\mathbb R^n)}ds.\\
\end{split}
\end{align}
For the last term on the right-hand side of \eqref{d30}, by \eqref{c12} we get
\begin{align}\label{d31}
\begin{split}
&6\mathbb E\left(\int^{t\wedge\tau_m}_{\tau}e^{2\eta s}\|u^{\varepsilon}(s)\|^2_{ H_{\rho}(\mathcal O)}\|\varpi_{\varepsilon}(u^{\varepsilon}(s),\mathcal L_{u^{\varepsilon}(s)})\|^2_{ L^2(l^2,H_{\rho}(\mathcal O))}ds\right)\\
&\le12\rho_2\mathbb E\left(\int^{t\wedge\tau_m}_{\tau}e^{2\eta s}\|u^{\varepsilon}(s)\|^2_{ H_{\rho}(\mathcal O)}\|\sigma_{1}\|^2_{L^2(\mathbb R^n,l^2)}ds\right)\\
&+
48\rho_2\|\kappa\|^2_{L^2(\mathbb R^n)}\|\beta\|^2_{l^2}\mathbb E\left(\int^{t\wedge\tau_m}_{\tau}e^{2\eta s}\|u^{\varepsilon}(s)\|^2_{ H_{\rho}(\mathcal O)}ds\right)\\
&+
48\|\kappa\|^2_{L^2(\mathbb R^n)}\|\beta\|^2_{l^2}\mathbb E\left(\int^{t\wedge\tau_m}_{\tau}e^{2\eta s}\|u^{\varepsilon}(s)\|^2_{ H_{\rho}(\mathcal O)}
\mathbb E(\|u^{\varepsilon}(s)\|^2_{ H_{\rho}(\mathcal O)})ds\right)\\
&+24\|\kappa\|^2_{L^{\infty}(\mathbb R^n)}\|\gamma\|^2_{l^2}\mathbb E\left(\int^{t\wedge\tau_m}_{\tau}e^{2\eta s}\|u^{\varepsilon}(s)\|^4_{ H_{\rho}(\mathcal O)}ds\right)\\
&\le\frac{1}{3}\eta\int^t_{\tau}e^{2\eta s}\mathbb E\left(\|u^{\varepsilon}(s)\|^4_{ H_{\rho}(\mathcal O)}\right)ds
+108\rho_2^2\frac{1}{\eta}\int^t_{\tau}e^{2\eta s}\|\sigma_{1}\|^4_{L^2(\mathbb R^n,l^2)}ds\\
&+48\rho_2\|\kappa\|^2_{L^2(\mathbb R^n)}\|\beta\|^2_{l^2}\int^{t}_{\tau}e^{2\eta s}\mathbb E\left(\|u^{\varepsilon}(s)\|^2_{ H_{\rho}(\mathcal O)}\right)ds\\
&+
48\|\kappa\|^2_{L^2(\mathbb R^n)}\|\beta\|^2_{l^2}\int^{t}_{\tau}e^{2\eta s}\mathbb E\left(\|u^{\varepsilon}(s)\|^4_{ H_{\rho}(\mathcal O)}\right)
ds\\
&+24\|\kappa\|^2_{L^{\infty}(\mathbb R^n)}\|\gamma\|^2_{l^2}\int^{t}_{\tau}e^{2\eta s}\mathbb E\left(\|u^{\varepsilon}(s)\|^4_{ H_{\rho}(\mathcal O)}\right)ds\\
&\le(\frac{1}{2}\eta+48\|\kappa\|^2_{L^2(\mathbb R^n)}\|\beta\|^2_{l^2}+24\|\kappa\|^2_{L^{\infty}(\mathbb R^n)}\|\gamma\|^2_{l^2})\int^{t}_{\tau}e^{2\eta s}\mathbb E\left(\|u^{\varepsilon}(s)\|^4_{ H_{\rho}(\mathcal O)}\right)ds\\
&+108\rho_2^2\frac{1}{\eta}\int^t_{\tau}e^{2\eta s}\|\sigma_{1}\|^4_{L^2(\mathbb R^n,l^2)}ds+1728\eta^{-2}\rho^2_2\|\kappa\|^4_{L^2(\mathbb R^n)}\|\beta\|^4_{l^2}e^{2\eta t}.
\end{split}
\end{align}
It follows from \eqref{d30}-\eqref{d31} that for all $t\ge\tau$,
\begin{align}\label{d32}
\begin{split}
&\mathbb E\left(e^{2\eta (t\wedge\tau_m)}\|u^{\varepsilon}(t\wedge\tau_m)\|^4_{ H_{\rho}(\mathcal O)}\right)+2(2\lambda-\eta)\mathbb E\left(\int^{t\wedge\tau_m}_{\tau}e^{2\eta s}\|u^{\varepsilon}(s)\|^4_{ H_{\rho}(\mathcal O)}ds\right)\\
&\le e^{2\eta \tau}\|\xi^{\varepsilon}\|^4_{ H_{\rho}(\mathcal O)}+12\frac{1}{\eta}\rho^2\int^{t}_{\tau}e^{2\eta s}\|\phi_1\|^2_{L^1(\mathbb R^n)}ds\\
&+108\rho_2^2\frac{1}{\eta}\int^t_{\tau}e^{2\eta s}\|\sigma_{1}\|^4_{L^2(\mathbb R^n,l^2)}ds+1728\eta^{-2}\rho^2_2\|\kappa\|^4_{L^2(\mathbb R^n)}\|\beta\|^4_{l^2}e^{2\eta t}\\
&+\int^{t}_{\tau}e^{2\eta s}\left(4\|\phi_1\|_{L^{\infty}(\mathbb R^n)}+4\|\psi_1\|_{L^1(\mathbb R^n)}+\frac{1}{3}\eta\right)\mathbb E\left(\|u^{\varepsilon}(s)\|^4_{H_{\rho}(\mathcal O)}\right)ds\\
&+(\frac{1}{2}\eta+48\|\kappa\|^2_{L^2(\mathbb R^n)}\|\beta\|^2_{l^2}+24\|\kappa\|^2_{L^{\infty}(\mathbb R^n)}\|\gamma\|^2_{l^2})\int^{t}_{\tau}e^{2\eta s}\mathbb E\left(\|u^{\varepsilon}(s)\|^4_{ H_{\rho}(\mathcal O)}\right)ds\\
&+2\eta\int^{t}_{\tau}e^{2\eta s}\mathbb E\left(\|u^{\varepsilon}(s)\|^4_{ H_{\rho}(\mathcal O)}\right)ds+\frac{27}{8\eta^3}\rho\int^{t}_{\tau}e^{2\eta s}\|g(s)\|^2_{L^2(\mathbb R^n)}ds.
\end{split}
\end{align}
Taking the limit of \eqref{d32} as $m\rightarrow\infty$, by Fatou's lemma we obtain for all $t\ge\tau$,
\begin{align}\label{d33}
\begin{split}
&\mathbb E\left(e^{2\eta t}\|u^{\varepsilon}(t)\|^4_{ H_{\rho}(\mathcal O)}\right)+2(2\lambda-\eta)\mathbb E\left(\int^{t}_{\tau}e^{2\eta s}\|u^{\varepsilon}(s)\|^4_{ H_{\rho}(\mathcal O)}ds\right)\\
&\le e^{2\eta \tau}\|\xi^{\varepsilon}\|^4_{ H_{\rho}(\mathcal O)}+12\frac{1}{\eta}\rho^2\int^{t}_{\tau}e^{2\eta s}\|\phi_1\|^2_{L^1(\mathbb R^n)}ds\\
&+108\rho_2^2\frac{1}{\eta}\int^t_{\tau}e^{2\eta s}\|\sigma_{1}\|^4_{L_2(\mathbb R^n,l^2)}ds+1728\eta^{-2}\rho^2_2\|\kappa\|^4_{L^2(\mathbb R^n)}\|\beta\|^4_{l^2}e^{2\eta t}\\
&+\int^{t}_{\tau}e^{2\eta s}\left(4\|\phi_1\|_{L^{\infty}(\mathbb R^n)}+4\|\psi_1\|_{L^1(\mathbb R^n)}+\frac{1}{3}\eta\right)\mathbb E\left(\|u^{\varepsilon}(s)\|^4_{H_{\rho}(\mathcal O)}\right)ds\\
&+(\frac{1}{2}\eta+48\|\kappa\|^2_{L^2(\mathbb R^n)}\|\beta\|^2_{l^2}+24\|\kappa\|^2_{L^{\infty}(\mathbb R^n)}\|\gamma\|^2_{l^2})\int^{t\wedge\tau_m}_{\tau}e^{2\eta s}\mathbb E\left(\|u^{\varepsilon}(s)\|^4_{ H_{\rho}(\mathcal O)}\right)ds\\
&+2\eta\int^{t}_{\tau}e^{2\eta s}\mathbb E\left(\|u^{\varepsilon}(s)\|^4_{ H_{\rho}(\mathcal O)}\right)ds+\frac{27}{8\eta^3}\rho\int^{t}_{\tau}e^{2\eta s}\|g(s)\|^2_{L^2(\mathbb R^n)}ds.
\end{split}
\end{align}
By \eqref{c38} and \eqref{d33} we get for all $t\ge\tau$,
\begin{align}\label{d34}
\begin{split}
&\mathbb E\left(\|u^{\varepsilon}(t)\|^4_{ H_{\rho}(\mathcal O)}\right)\\
&\le e^{2\eta (\tau-t)}\|\xi^{\varepsilon}\|^4_{ H_{\rho}(\mathcal O)}+12\frac{1}{\eta}\rho^2\int^{t}_{\tau}e^{2\eta (s-t)}\|\phi_1\|^2_{L^1(\mathbb R^n)}ds\\
&+108\rho_2^2\frac{1}{\eta}\int^t_{\tau}e^{2\eta (s-t)}\|\sigma_{1}\|^4_{L^2(\mathbb R^n,l^2)}ds+1728\eta^{-2}\rho^2_2\|\kappa\|^4_{L^2(\mathbb R^n)}\|\beta\|^4_{l^2}\\
&+\frac{27}{8\eta^3}\rho\int^{t}_{\tau}e^{2\eta (s-t)}\|g(s)\|^2_{L^2(\mathbb R^n)}ds\\
&\le e^{2\eta (\tau-t)}\|\xi^{\varepsilon}\|^4_{ H_{\rho}(\mathcal O)}+6\frac{1}{\eta^2}\rho^2\|\phi_1\|^2_{L^1(\mathbb R^n)}+
54\rho_2^2\frac{1}{\eta^2}\|\sigma_{1}\|^4_{L^2(\mathbb R^n,l^2)}\\
&+1728\eta^{-2}\rho^2_2\|\kappa\|^4_{L^2(\mathbb R^n)}\|\beta\|^4_{l^2}
+\frac{27}{16\eta^4}\rho\|g_0\|^2_{C^b(\mathbb R,L^2(\mathbb R^n))}.
\end{split}
\end{align}
Since $\mathbb E(\|\xi^{\varepsilon}\|^4_{ H_{\rho}(\mathcal O)})\le R$, we have
$$\lim\limits_{t\rightarrow\infty}e^{2\eta (\tau-t)}\|\xi^{\varepsilon}\|^4_{ H_{\rho}(\mathcal O)}\le\lim\limits_{t\rightarrow\infty}e^{2\eta (\tau-t)}R=0,$$
and hence there exists $T=T(R)$ implies that for all $t-\tau>T$,
$$e^{2\eta (\tau-t)}\|\xi^{\varepsilon}\|^4_{ H_{\rho}(\mathcal O)}\le 1,$$
which along with \eqref{d34} implies that for all $t-\tau>T$.
\begin{align}
\begin{split}
&\mathbb E\left(\|u^{\varepsilon}(t)\|^4_{ H_{\rho}(\mathcal O)}\right)\le1+6\frac{1}{\eta^2}\rho^2\|\phi_1\|^2_{L^1(\mathbb R^n)}\\
&+54\rho_2^2\frac{1}{\eta^2}\|\sigma_{1}\|^4_{L^2(\mathbb R^n,l^2)}+1728\eta^{-2}\rho^2_2\|\kappa\|^4_{L^2(\mathbb R^n)}\|\beta\|^4_{l^2}\\
&+\frac{27}{16\eta^4}\rho\|g_0\|^2_{C^b(\mathbb R,L^2(\mathbb R^n))},
\end{split}
\end{align}
which completes the proof.
\end{proof}
\end{lm}
\section{Existence of Uniform Measure Attractors}

\setcounter{equation}{0}
In the section, we prove the existence and uniqueness of uniform measure attractor of \eqref{c8} in $\mathcal P_4(L^2(\mathcal O))$. Firstly, we define a process
in $\mathcal P_4(L^2(\mathcal O))$.

Given $t \geq r$, for every $\mu \in \mathcal P_4(L^2(\mathcal O))$, define $P^{g,\varepsilon}_*:\mathcal P_4(L^2(\mathcal O))\rightarrow\mathcal P_4(L^2(\mathcal O))$ by
\begin{align}\label{e1}
P_*^{g,\varepsilon}(t, \tau) \mu=\mathcal L_{u^{g,\varepsilon}(t,\tau,\xi^{\varepsilon})},
\end{align}
where $u^{g,\varepsilon}(t,\tau,\xi^{\varepsilon})$ is the solution of \eqref{c8} with $\xi^{\varepsilon}\in L^4(\Omega,\mathcal F_{\tau},L^2(\mathcal O))$ such that
$\mathcal L_{\xi^{\varepsilon}}=\mu$. In terms of \eqref{e1}, for every $t\in\mathbb R^+$ and $\tau\in\mathbb R$, define $U^{g,\varepsilon}(t,\tau):\mathcal P_4(L^2(\mathcal O))\rightarrow\mathcal P_4(L^2(\mathcal O))$ by, for every $\mu\in\mathcal P_4(L^2(\mathcal O))$,
\begin{align}\label{e18}
U^{g,\varepsilon}(\tau+t,\tau)\mu=P_*^{g,\varepsilon}(\tau+t,\tau)\mu.
\end{align}
By the uniqueness of solutions for \eqref{c8}, the operator $U^{g,\varepsilon}(t,\tau)$ satisfies the multiplicative properties:
$$U^{g,\varepsilon}(t,\tau)=U^{g,\varepsilon}(t,s)U^{g,\varepsilon}(s,\tau).$$
for all $t\ge s\ge\tau,$ $\tau\in\mathbb R.$
$$U^{g,\varepsilon}(\tau,\tau)=I,\quad \tau\in\mathbb R,$$
where $I$ is the identity operator. Furthermore, the following translation identity holds by a similar argument to that of Lemma 4.1 in \cite{LWW2021}
$$U^{g,\varepsilon}(t+h,\tau+h)=U^{T(h)g,\varepsilon}(t,\tau),$$
for all $h\in\mathbb R$, $t\ge\tau,$ $\tau\in\mathbb R$.

In a manner analogous to the preceding, we may likewise define a process designated as $U^{g,0}$ in accordance with \eqref{c9}.

Next, we establish the continuity of  $U^{g,\varepsilon}(t,\tau)$ with respect to the topology of $\mathcal P_4(L^2(\mathcal O))\times \mathcal H(g_0)$.
\begin{lm}\label{e17}
Suppose $\mathbf(A1)$-$\mathbf(A3)$ hold. Let $\xi^{\varepsilon}, \xi^{\varepsilon}_n\in L^4(\Omega,\mathcal F_{\tau},L^2(\mathcal O))$ such that
$\mathbb E\left(\|\xi^{\varepsilon}\|^4\right)\le R$ and $\mathbb E\left(\|\xi^{\varepsilon}_n\|^4\right)\le R$ for some $R>0$.
If $\mathcal L_{\xi^{\varepsilon}_n}\rightarrow\mathcal L_{\xi^{\varepsilon}}$ weakly and $g_n\rightarrow g$ in $\mathcal H(g_0)$, then
 for every $\tau\in\mathbb R$, $t\ge\tau$ and $0<\varepsilon<\varepsilon_0$,
$\mathcal L_{u^{g_n,\varepsilon}(t,\tau,\xi^{\varepsilon}_n)}\rightarrow\mathcal L_{u^{g,\varepsilon}(t,\tau,\xi^{\varepsilon})}$ weakly.
\begin{proof}
Since $\mathcal L_{\xi^{\varepsilon}_n}\rightarrow\mathcal L_{\xi^{\varepsilon}}$ weakly, by the Skorokhov theorem, there exist a probability space $(\widetilde{\Omega}, \widetilde{\mathcal{F}}, \widetilde{\mathbb{P}})$ and random variables $\widetilde{\xi^{\varepsilon}}$ and $\widetilde{\xi^{\varepsilon}_{n}}$ defined in $(\widetilde{\Omega}, \widetilde{\mathcal{F}}, \widetilde{\mathbb{P}})$ such that the distributions of $\widetilde{\xi^{\varepsilon}}$ and $\widetilde{\xi^{\varepsilon}_{n}}$ coincide with that of $\xi^{\varepsilon}$ and $\xi^{\varepsilon}_n$, respectively.
Furthermore, $\widetilde{\xi^{\varepsilon}}_{n} \rightarrow \widetilde{\xi^{\varepsilon}}$ $\widetilde{\mathbb{P}}$-almost surely. Note that $\widetilde{\xi^{\varepsilon}}$, $\widetilde{\xi^{\varepsilon}_{n}}$ and $W$ can be considered as random variables defined in the product space $(\Omega \times \widetilde{\Omega}, \mathcal{F} \times \widetilde{\mathbb{F}}, \mathbb{P} \times \widetilde{\mathbb{P}})$. So we may consider the solutions of the stochastic equation in the product space with initial data $\widetilde{\xi^{\varepsilon}}$ and $\widetilde{\xi^{\varepsilon}_{n}}$, instead of the solutions in $(\Omega, \mathcal{F}, \mathbb{P})$ with initial data $\xi^{\varepsilon}$ and $\xi^{\varepsilon}_n$. However, for simplicity, we will not distinguish the new random variables from the original ones, and just consider the solutions of the equation in the original space. Since $\widetilde{\xi^{\varepsilon}_{n}} \rightarrow \widetilde{\xi^{\varepsilon}}$ $(\mathbb{P} \times \widetilde{\mathbb{P}})$-almost surely, without loss of generality, we simply assume that $\xi^{\varepsilon}_n \rightarrow \xi^{\varepsilon}$ $\mathbb{P}$-almost surely.

Let $u^{g_n,\varepsilon}_n(t, \tau)=u^{g_n,\varepsilon}\left(t, \tau, \xi^{\varepsilon}_n \right)$, $u^{g,\varepsilon}(t, \tau)=u^{g,\varepsilon}\left(t, \tau, \xi^{\varepsilon}\right)$ and $\\ v_n^{\varepsilon}(t, \tau)=u^{g_n,\varepsilon}\left(t, \tau, \xi^{\varepsilon}_n \right)-u^{g,\varepsilon}\left(t, \tau, \xi^{\varepsilon}\right)$. Then by \eqref{c8} we have, for all $t \geq \tau$,

$$
\begin{gathered}
d v^{\varepsilon}_n(t)-A_{\varepsilon} v^{\varepsilon}_n(t) d t+\lambda v^{\varepsilon}_n(t) d t+\left(f_{\varepsilon}\left(y, u^{g_n,\varepsilon}_n(t), \mathcal{L}_{u^{g_n,\varepsilon}_n(t)}\right)-f_{\varepsilon}\left(y, u^{g,\varepsilon}(t), \mathcal{L}_{u^{g,\varepsilon}(t)}\right)\right) d t \\
=\left(g_n\left(t, y^*\right)-g\left(t, y^*\right) \right)d t+\left(\varpi_{\varepsilon}\left( u^{g_n,\varepsilon}_n(t), \mathcal{L}_{u^{g_n,\varepsilon}_n(t)}\right)-\varpi_{\varepsilon}\left( u^{g,\varepsilon}(t), \mathcal{L}_{u^{g,\varepsilon}(t)}\right)\right) d W(t).
\end{gathered}
$$

By Ito's formula  we have for all $t\ge\tau$,
\begin{align}\label{e13}
\begin{split}
&\left\|v^{\varepsilon}_n(t)\right\|^2_{H_{\rho}(\mathcal O)}+2 \int_\tau^t a_{\varepsilon}(v^{\varepsilon}_n(s),v^{\varepsilon}_n(s))ds+2 \lambda \int_\tau^t\left\|v^{\varepsilon}_n(s)\right\|^2_{H_{\rho}(\mathcal O)} d s\\
&+2\int_\tau^t\int_{\mathcal O}\rho \left(f_{\varepsilon}\left(y, u^{g_n,\varepsilon}_n(s), \mathcal{L}_{u^{g_n,\varepsilon}_n(s)}\right)-f_{\varepsilon}\left( y, u^{g,\varepsilon}(s), \mathcal{L}_{u^{g,\varepsilon}(s)}\right)\right) v^{\varepsilon}_n(s)dyds\\
&=\|\xi^{\varepsilon}_n-\xi^{\varepsilon}\|^2_{H_{\rho}(\mathcal O)}+2\int^t_{\tau}\left(g_n\left(s, \cdot\right)-g\left(s, \cdot\right),v_n^{\varepsilon}(s) \right)_{H_{\rho}(\mathcal O)}ds\\
&+\int^t_{\tau}\|\varpi_{\varepsilon}\left( u^{g_n,\varepsilon}_n(s), \mathcal{L}_{u^{g_n,\varepsilon}_n(s)}\right)-\varpi_{\varepsilon}\left( u^{g,\varepsilon}(s), \mathcal{L}_{u^{g,\varepsilon}(s)}\right)\|^2_{L_2(l^2,H_{\rho}(\mathcal O))}ds\\
&+2\int^t_{\tau}\left(v^{\varepsilon}_n(s),\varpi_{\varepsilon}\left( u^{g_n,\varepsilon}_n(s), \mathcal{L}_{u^{g_n,\varepsilon}_n(s)}\right)-\varpi_{\varepsilon}\left( u^{g,\varepsilon}(s), \mathcal{L}_{u^{g,\varepsilon}(s)}\right)\right)dW.
\end{split}
\end{align}
By \eqref{c4} and \eqref{c14} we have
\begin{align}
\begin{split}
&-2\int_\tau^t\int_{\mathcal O}\rho \left(f_{\varepsilon}\left(y, u^{g_n,\varepsilon}_n(s), \mathcal{L}_{u^{g_n,\varepsilon}_n(s)}\right)-f_{\varepsilon}\left(y, u^{g,\varepsilon}(t), \mathcal{L}_{u^{g,\varepsilon}(t)}\right)\right) v^{\varepsilon}_n(s)dyds\\
&=-2\int_\tau^t\int_{\mathcal O}\rho \left(f_{\varepsilon}\left( y,u^{g_n,\varepsilon}_n(s), \mathcal{L}_{u^{g_n,\varepsilon}_n(s)}\right)-f_{\varepsilon}\left(y, u^{g,\varepsilon}(s), \mathcal{L}_{u^{g_n,\varepsilon}_n(s)}\right)\right)v^{\varepsilon}_ndyds\\
&-2\int_\tau^t\int_{\mathcal O}\rho \left(f_{\varepsilon}\left(y, u^{g,\varepsilon}(s), \mathcal{L}_{u^{g_n,\varepsilon}_n(s)}\right)-f_{\varepsilon}\left( y, u^{g,\varepsilon}(t), \mathcal{L}_{u^{g,\varepsilon}(t)}\right)\right)v^{\varepsilon}_n(s)dyds\\
&\le\int_\tau^t\left((2\|\phi_4\|_{L^{\infty}(\mathbb R^n)}+\|\phi_3\|_{L^{\infty}(\mathbb R^n)})\|v^{\varepsilon}_n(s)\|_{H_{\rho}(\mathcal O)}^2
+\mathbb E(\|v^{\varepsilon}_n(s)\|_{H_{\rho}(\mathcal O)}^2)\|\phi_3\|_{L^{1}(\mathbb R^n)}\right)ds.
\end{split}
\end{align}
By Young inequality, we have
\begin{align}
\begin{split}
&2\int^t_{\tau}\left(g_n\left(s, \cdot\right)-g\left(s, \cdot\right),v_n^{\varepsilon}(s) \right)_{H_{\rho}(\mathcal O)}ds\\
&\le \int^t_{\tau}\|g_n-g\|^2_{C_b(\mathbb R,L^2(\mathbb R^n)}+\|v_n^{\varepsilon}(s)\|^2_{H_{\rho}(\mathcal O)}ds\\
&\le \|g_n-g\|^2_{C_b(\mathbb R,L^2(\mathbb R^n)}(t-\tau)+\int^t_{\tau}\|v_n^{\varepsilon}(s)\|^2_{H_{\rho}(\mathcal O)}ds.
\end{split}
\end{align}
By \eqref{c18} we have
\begin{align}\label{e14}
\begin{split}
&\int^t_{\tau}\|\varpi_{\varepsilon}\left( u^{g_n,\varepsilon}_n(s), \mathcal{L}_{u^{g_n,\varepsilon}_n(s)}\right)-\varpi_{\varepsilon}\left( u^{g,\varepsilon}(s), \mathcal{L}_{u^{g,\varepsilon}(s)}\right)\|^2_{L^2(l^2,H_{\rho}(\mathcal O))}ds\\
&\le
2\|L_{\varpi}\|^2_{l^2}\int^t_{\tau}\left(\|\kappa\|^2_{L^{\infty}(\mathbb R^n)}\|v_n^{\varepsilon}(s)\|^2_{H_{\rho}(\mathcal O)}
+\|\kappa\|^2_{L^2(\mathbb R^n)}\mathbb \|v_n^{\varepsilon}(s)\|^2_{H_{\rho}(\mathcal O)}\right)ds.
\end{split}
\end{align}
By \eqref{e13}-\eqref{e14}, we find  that for every $T>0$, such that for all $t\in[\tau,\tau+T]$,
\begin{align}\label{e15}
\begin{split}
&\mathbb E\left(\|v^{\varepsilon}_n(t)\|^2_{H_{\rho}(\mathcal O)}\right)\le \mathbb E\left(\|\xi^{\varepsilon}_n-\xi^{\varepsilon}\|^2_{H_{\rho}(\mathcal O)}\right)+\|g_n(s)-g(s)\|^2_{C_b(\mathbb R,L^2(\mathbb R^n)}T\\
&+(1+2\|\phi_4\|_{L^{\infty}(\mathbb R^n)}+\|\phi_3\|_{L^{\infty}(\mathbb R^n)}+\|\phi_3\|_{L^{1}(\mathbb R^n)})\int_\tau^t\mathbb E\left(\|v^{\varepsilon}_n(s)\|_{H_{\rho}(\mathcal O)}^2\right)ds\\
&+2\|L_{\varpi}\|^2_{l^2}(\|\kappa\|^2_{L^{\infty}(\mathbb R^n)}+\|\kappa\|^2_{L^2(\mathbb R^n)})\int^t_{\tau}\mathbb E\left(
\|v_n^{\varepsilon}(s)\|^2_{H_{\rho}(\mathcal O)}\right)ds.
\end{split}
\end{align}
By \eqref{e15} and Gronwall's lemma, we obtain, for all $t\in[\tau,\tau+T]$,
\begin{align}\label{e16}
\begin{split}
&\mathbb E\left(\|v^{\varepsilon}_n(t)\|^2_{H_{\rho}(\mathcal O)}\right)\le\left(\mathbb E\left(\|\xi^{\varepsilon}_n-\xi^{\varepsilon}\|^2_{H_{\rho}(\mathcal O)}\right)+\|g_n-g\|^2_{C_b(\mathbb R,L^2(\mathbb R^n))}T\right)e^{c_1(t-\tau)},
\end{split}
\end{align}
where $c_1>0$ is a constant independent of $n,\tau$ and $t$.
Since $\mathbb E(\|\xi^{\varepsilon}_n\|^4_{H_{\rho}(\mathcal O)})\le R$, we see that the sequence $\{\xi^{\varepsilon}_n\}^{\infty}_{n=1}$ is uniformly integrable in $L^2(\Omega, L^2(\mathcal O))$. Then using the assumption that $\xi^{\varepsilon}_n\rightarrow\xi^{\varepsilon}$ $\mathbb P$-almost surely, we obtain from Vitali's
theorem that $\xi^{\varepsilon}_n\rightarrow\xi^{\varepsilon}$ in $L^2(\Omega,L^2(\mathcal O))$ and  $g_n\rightarrow g$ in $\mathcal H(g_0)$, which along with \eqref{e16} shows that
$u^{g_n,\varepsilon}\left(t, \tau, \xi^{\varepsilon}_n \right)\rightarrow u^{g,\varepsilon}\left(t, \tau, \xi^{\varepsilon}\right)$ in $L^2(\Omega,L^2(\mathcal O))$
and  hence also in distribution.
\end{proof}
\end{lm}

By Lemma \ref{e17}, we find that the process $U^{g,\varepsilon}$ given by \eqref{e18} is continuous  over bounded of $\mathcal P_4(L^2(\mathcal O))\times \mathcal H(g_0)$.
\begin{lm}\label{e32}
Suppose $\mathbf(A1)$-$\mathbf(A3)$ and \eqref{c10} hold. Denoted by
\begin{equation}\label{d9}
K=B_{\mathcal P_4(L^2(\mathcal O))}(G^{\frac{1}{4}}_1),
\end{equation}
where
\begin{align}
G_1=\frac{M_7}{\rho_1},
\end{align}
$M_7>0$ is the same constants as in Lemma \ref{d35}. Then $K$ is a closed uniform absorbing set of the family of  $\{U^{g_{},\varepsilon}\}_{g_{}\in\mathcal H( g_0)}$.
\begin{proof}
First, note that $K$ is a closed subset of $\mathcal P_4(L^2(\mathcal O))$. Subsequently,
by \eqref{d9} and Lemma \ref{d35}, we see that then for every $R>0$, there exists $T=T(R)>0$ such that
$\{U^{g_{},\varepsilon}(t,0)\}_{g_{}\in\mathcal H(g_0)}$ satisfies
$$U^{g_{},\varepsilon}(t,0)B_{\mathcal P_4(L^2(\mathcal O))}(R)\subseteq K,\ \text{for all}\ g_{}\in\mathcal H(g_0)\ \text{and}
\ t\ge T.$$

\end{proof}
\end{lm}
We now present the uniformly asymptotically compact of the family of process $\{U^{g_{},\varepsilon}\}_{g_{}\in\mathcal H(g_0)}$
with respect to $g_{}\in\mathcal H(g_0)$.
\begin{lm}\label{e2}
Suppose $\mathbf(A1)$-$\mathbf(A3)$ and \eqref{c10} hold. Then the family of processes $\\ \{U^{g_{},\varepsilon}(t,0)\}_{g_{}\in\mathcal H(g_0)}$ is
uniformly asymptotically compact in $\mathcal P_4(L^2(\mathcal O))$; that is, $\\ \{U^{g_{n},\varepsilon}(t_n,0)\mu_n\}$ has a convergent subsequence in
$\mathcal P_4(L^2(\mathcal O))$ whenever $t_n\rightarrow+\infty$ and $(\mu_n,g_{n})$ is bounded in $\mathcal P_4(L^2(\mathcal O))\times \mathcal H(g_0)$.
\begin{proof}
Given $v_n\in L^4(\Omega,\mathcal F_0,L^2(\mathcal O))$ with distribution $\mu_n$, i.e. $\mathcal L_{v_n}=\mu_n$, we consider the solution $u^{g_n,\varepsilon}(t,0,v_n)$ of equation \eqref{c8}
with initial data $v_n$ at initial time $0$.
To complete the proof, by Prohorov theorem, it is to prove that the sequence $ \{\mathcal L_{(u^{g_n,\varepsilon}(t_n,0,v_n)}\}$ is tight in $L^2(\mathcal O)$.

Let $\theta$ be a cut-off smooth function given by \eqref{d36}, $\theta_m=\theta(\frac{|y^*|^2}{m^2})$ for every $m\in\mathbb N$ and $y^*\in\mathbb R^n$.
Then the solution $u^{g,\varepsilon}$ can be decomposed as:
$$u^{g_n,\varepsilon}(t_n,0,v_n)=\theta_m u^{g_n,\varepsilon}(t_n,0,v_n)+(1-\theta_m) u^{g_n,\varepsilon}(t_n,0,v_n).$$
Note that there exists $c_1=c_1(\theta)>0$ such that
$$\sup_{s\in\mathbb R^+}|\theta'(s)|\le c_1,$$
and hence for all $m,n\in\mathbb N$,
\begin{align}\label{e26}
\|(1-\theta_m) u^{g_n,\varepsilon}(t_n,0,v_n)\|^2_{H^1_{\varepsilon}(\mathcal O)}\le c_2\| u^{g_n,\varepsilon}(t_n,0,v_n)\|^2_{H^1_{\varepsilon}(\mathcal O)},
\end{align}
where $c_2=1+c_1^2$.

By Lemma \ref{d10}, we know that for every $0<\varepsilon<\varepsilon_0$, there exists $N_1\in\mathbb N$, independent of $\varepsilon$, and $c_3=c_3(\varepsilon,g_0)$ such that for all
$n\ge N_1$,
\begin{align}\label{e27}
\mathbb E(\|u^{g_{n},\varepsilon}(t_n,0,v_n)\|^2_{H^{1}_{\varepsilon}(\mathcal O)})\le c_3.
\end{align}
By \eqref{e26} and \eqref{e27}, we have for all $m\in\mathbb N$ and $n\ge N_1$,
\begin{align}\label{e28}
\mathbb E\left(\|(1-\theta_m) u^{g_n,\varepsilon}(t_n,0,v_n)\|^2_{H^1_{\varepsilon}(\mathcal O)}\right)\le c_2c_3.
\end{align}
By \eqref{e28} and Chebyshev's inequality, we get for all $m\in\mathbb N$ and $n\ge N_1$,
\begin{align}
\mathbb {P}(\|(1-\theta_m) u^{g_n,\varepsilon}(t_n,0,v_n)\|^2_{H^1_{\varepsilon}(\mathcal O)}\ge R_1)\le \frac{c_2c_3}{R^2_1}\rightarrow0\ \text{as}\ R_1\rightarrow\infty.
\end{align}
Therefore, for every $\delta>0$, there exists $R_2=R_2(\delta)>0$ such that for all $m\in\mathbb N$ and $n\ge N_1$,
\begin{align}\label{e29}
\mathbb P(\|(1-\theta_m) u^{g_n,\varepsilon}(t_n,0,v_n)\|^2_{H^1_{\varepsilon}(\mathcal O)}> R_2)<\delta.
\end{align}
Let
$$Z_{\delta}=\{u^{g_n,\varepsilon}\in H^1_{\varepsilon}(\mathcal O):\|u^{g_n,\varepsilon}\|_{H^1_{\varepsilon}(\mathcal O)}\le R_2(\delta);\quad u^{g_n,\varepsilon}(y)=0\
\text{for a.e.}\ |y^*|>\sqrt{2}m\}.$$
Then $Z_{\delta}$ is a compact subset of $L^2(\mathcal O)$. By \eqref{e29} we have for all $m\in\mathbb N$ and $n\ge N_1$,
\begin{align}\label{e30}
\mathbb P(\{(1-\theta_m) u^{g_n,\varepsilon}(t_n,0,v_n)\in Z_{\delta}\})>1-\delta.
\end{align}
Since $\delta>0$ is arbitrary, by \eqref{e30} we find that for every $m\in\mathbb N$,
\begin{align}\label{e31}
\{\mathcal L_{(1-\theta_m) u^{g_n,\varepsilon}(t_n,0,v_n)}\}^{\infty}_{n=1} \text{ is tight in}\  L^2(\mathcal O),
\end{align}
where $\mathcal L_{(1-\theta_m) u^{g_n,\varepsilon}(t_n,0,v_n)}$ is the distribution of $(1-\theta_m) u^{g_n,\varepsilon}(t_n,0,v_n)$ in $ L^2(\mathcal O)$.

Next, we demonstrate that the sequence $\{\mathcal L_{u^{g_n,\varepsilon}(t_n,0,v_n))}\}$ is tight in $L^2(\mathcal O)$ by utilizing uniform tail-estimates. Indeed,
by invoking Lemma \ref{d13}, we deduce that for every $\delta_1>0$, there exists  $N_2=N_2(\delta)\in\mathbb N$ and $m_0=m_0(\delta)\in\mathbb N$  such that for all $n\ge N_2$,
\begin{align}\label{e19}
\mathbb E\left(\int_{|y^*|\ge m_0}\int^1_0|u^{g_{n},\varepsilon}(t_n,0,v_n)(y^*,y_{n+1})|^2d{y_{n+1}}dy^*)\right)\le\frac{1}{9}\delta_1^2.
\end{align}
By \eqref{e19} we get, for all $n\ge N_2$,
\begin{align}\label{e22}
\mathbb E\left(\|\theta_{m_0} u^{g_{n},\varepsilon}(t_n,0,v_n)\|^2\right)<\frac{1}{9}\delta^2_1.
\end{align}
By \eqref{e22}, we see that $\{\mathcal L_{(1-\theta_{m_0}) u^{g_n,\varepsilon}(t_n,0,v_n)}\}^{\infty}_{n=N_2} \text{ is tight in}\  L^2(\mathcal O)$, and hence there exists
$n_1,\cdots,n_l\ge N_2$ such that
\begin{align}\label{e20}
\{\mathcal L_{(1-\theta_{m_0}) u^{g_n,\varepsilon}(t_n,0,v_n)}\}^{\infty}_{n=N_2} \subset
\bigcup\limits^l_{j=1}B\left(\mathcal L_{(1-\theta_{m_0}) u^{g_n,\varepsilon}(t_{n_j},0,v_{n_j})},\frac{1}{3}\delta_1)\right),
\end{align}
where $B\left(\mathcal L_{(1-\theta_{m_0}) u^{g_n,\varepsilon}(t_{n_j},0,v_{n_j})},\frac{1}{3}\delta_1\right)$ is the $\frac{1}{3}\delta_1$-neighborhood of
$\\ \mathcal L_{(1-\theta_{m_0}) u^{g_n,\varepsilon}(t_{n_j},0,v_{n_j})}$ in the space $\mathcal P_4(L^2(\mathcal O),d_{\mathcal P(L^2(\mathcal O))})$.
We claim:
\begin{align}\label{e23}
\{\mathcal L_{ u^{g_n,\varepsilon}(t_n,0,v_n)}\}^{\infty}_{n=N_2}\subset\bigcup\limits^l_{j=1}B\left(\mathcal L_{ u^{g_n,\varepsilon}(t_{n_j},0,v_{n_j})},\delta_1\right).
\end{align}
Given $n\ge N_2$ by \eqref{e20} we know that there exist $j\in\{1,\cdots,l\}$ such that
\begin{align}\label{e21}
\mathcal L_{(1-\theta_{m_0}) u^{g_n,\varepsilon}(t_n,0,v_n)}\in B\left(\mathcal L_{(1-\theta_{m_0}) u^{g_n,\varepsilon}(t_{n_j},0,v_{n_j})},\frac{1}{3}\delta_1\right).
\end{align}
By \eqref{e22} and \eqref{e21} we have
\begin{align}
\begin{split}
&d_{\mathcal P(L^2(\mathcal O))}\left(\mathcal L_{u^{g_{n},\varepsilon}(t_n,0,v_n)},\mathcal L_{ u^{g_n,\varepsilon}(t_{n_j},0,v_{n_j})}\right)\\
&=\sup_{\substack{\varphi\in L_b(X)\\ \|\varphi\|_L\le 1}}\left|\int_{L^2(\mathcal O)}\varphi d\mathcal L(u^{g_{n},\varepsilon}(t_n,0,v_n)-
\int_{L^2(\mathcal O)}\varphi d\mathcal L( u^{g_n,\varepsilon}(t_{n_j},0,v_{n_j})\right|\\
&\le\sup_{\substack{\varphi\in L_b(X)\\ \|\varphi\|_L\le 1}}\left|\mathbb E(\varphi(u^{g_{n},\varepsilon}(t_n,0,v_n)))-\mathbb E(\varphi( u^{g_n,\varepsilon}(t_{n_j},0,v_{n_j})))\right|\\
&\le\sup_{\substack{\varphi\in L_b(X)\\ \|\varphi\|_L\le 1}}\left|\mathbb E(\varphi(u^{g_{n},\varepsilon}(t_n,0,v_n)))-
\mathbb E(\varphi(1-\theta_{m_0})(u^{g_{n},\varepsilon}(t_n,0,v_n)))\right|\\
&+\sup_{\substack{\varphi\in L_b(X)\\ \|\varphi\|_L\le 1}}\left|\mathbb E(\varphi(1-\theta_{m_0})(u^{g_{n},\varepsilon}(t_n,0,v_n)))-
\mathbb E(\varphi(1-\theta_{m_0})(u^{g_{n},\varepsilon}(t_{n_j},0,v_{n_j})))\right|\\
&+\sup_{\substack{\varphi\in L_b(X)\\ \|\varphi\|_L\le 1}}\left|\mathbb E(\varphi(1-\theta_{m_0})(u^{g_{n},\varepsilon}(t_{n_j},0,v_{n_j})))-
\mathbb E(\varphi(u^{g_{n},\varepsilon}(t_{n_j},0,v_{n_j})))\right|\\
&\le\mathbb E(\|\theta_{m_0}u^{g_{n},\varepsilon}(t_n,0,v_n)\|)+\mathbb E(\|\theta_{m_0}(u^{g_{n},\varepsilon}(t_{n_j},0,v_{n_j})\|)\\
&+d_{\mathcal P(L^2(\mathcal O))}\left(\mathcal L((1-\theta_{m_0})u^{g_{n},\varepsilon}(t_n,0,v_n),\mathcal L( (1-\theta_{m_0})u^{g_n,\varepsilon}(t_{n_j},0,v_{n_j})\right)\\
&\le\frac{1}{3}\delta_1+\frac{1}{3}\delta_1+\frac{1}{3}\delta_1=\delta_1,
\end{split}
\end{align}
which leads to the conclusion stated in \eqref{e23}. Since $\delta_1>0$ is arbitrary, we infer from \eqref{e23} that the sequence $\mathcal L_{u^{g_{n},\varepsilon}(t_n,0,v_n)}$ is tight in $L^2(\mathcal O)$. Consequently, there exists $\nu\in\mathcal P(L^2(\mathcal O))$ such that, possibly along a subsequence,
\begin{align}\label{e24}
\mathcal L_{u^{g_{n},\varepsilon}(t_n,0,v_n)}\rightarrow \nu \ \text{weakly}.
\end{align}
It remains to show $\nu\in\mathcal P_4(L^2(\mathcal O))$. Let $K$ be the closed uniform absorbing set of  $\{U^{g_{},\varepsilon}\}_{g_{}\in\mathcal H( g_0)}$ given by
\eqref{d9}. Then there exists $N_3\in\mathbb N$ such that for all $n\ge N_3$
\begin{align}\label{e25}
\mathcal L_{u^{g_{n},\varepsilon}(t_n,0,v_n)}\in K.
\end{align}
Since $K$ is closed with respect to the weak topology of $\mathcal P_4(L^2(\mathcal O))$, by \eqref{e24}-\eqref{e25} we obtain $\nu\in K$ and
thus $\nu\in\mathcal P_4(L^2(\mathcal O))$. This completes the proof.
\end{proof}
\end{lm}
\begin{thm}\label{e3}
Suppose $\mathbf(A1)$-$\mathbf(A3)$ and \eqref{c10} hold. Then the family of processes $\\ \{U^{g,\varepsilon}(t,\tau)\}_{g_{}\in\mathcal H(g_0)}$ associated with
\eqref{c8} has a unique uniform measure attractor $\mathcal A_{\varepsilon}$ in $\mathcal P_4(L^2(\mathcal O))$, which is given by
$$\mathcal A_{\varepsilon}=\underset{g \in \mathcal{H}\left(g_0\right)}{\cup} \mathcal{K}_{g,\varepsilon}(0).
$$
\begin{proof}
According to Lemma  \ref{e17}, the family of processes $\{U^{g,\varepsilon}(t,\tau)\}_{g_{}\in\mathcal H(g_0)}$ possesses joint continuity
over bounded of $\mathcal P_4(L^2(\mathcal O))$ and $\mathcal H(g_0)$.
Leveraging Lemmas \ref{e17}, \ref{e32} and \ref{e2}, Theorem \ref{b1} subsequently yields the existence and uniqueness of the uniform measure attractors for the aforementioned family of processes $\{U^{g,\varepsilon}(t,\tau)\}_{g_{}\in\mathcal H(g_0)}$.
\end{proof}
\end{thm}
The next result is concerned with existence and uniqueness of  uniform measure attractors
for $\{U^{g,0}(t,\tau)\}_{g_{}\in\mathcal H(g_0)}$ associated with problem \eqref{c9} which is analogous to Theorem \ref{e3}.
\begin{thm}
Suppose $\mathbf(A1)$-$\mathbf(A3)$ and \eqref{c10} hold. Then the family of processes $\\ \{U^{g,0}(t,\tau)\}_{g_{}\in\mathcal H(g_0)}$ associated with
\eqref{c9} has a unique uniform measure attractor $\mathcal A_{0}$ in $\mathcal P_4(L^2(\mathbb R^n))$ which is given by
$$\mathcal A_{0}=\underset{g \in \mathcal{H}\left(g_0\right)}{\cup} \mathcal{K}_{g,0}(0).
$$
\end{thm}

\section{Upper semicontinuity of uniform measure attractors}
\setcounter{equation}{0}
In this section, we prove the upper semicontinuity of uniform measure attractors for the non-autonomous stochastic reaction-diffusion equations when the $(n+1)$-dimensional thin domains collapse to an $n$-dimensional domain. To that end, we need the average operator $\mathcal{M}: L^2(\mathcal{O}) \rightarrow$ $L^2(\mathbb R^n)$ as given by: for every $\varphi \in L^2(\mathcal{O})$,
$$
\mathcal{M} \varphi\left(y^*\right)=\int_0^1 \varphi\left(y^*, y_{n+1}\right) d y_{n+1}, \quad y^* \in \mathbb R^n .
$$
Let $\mathcal{I}: L^2(\mathbb R^n) \rightarrow L^2(\mathcal{O})$ be the operator given by: for every $\varphi \in L^2(\mathbb R^n)$,
$$
\mathcal{I} \varphi(y)=\varphi\left(y^*\right), \quad \forall y=\left(y^*, y_{n+1}\right) \in \mathcal{O} .
$$
The following property of the operator $\mathcal{M}$ from \cite{HR01} will be used in the sequel.
\begin{lm}\label{f6}
If $u \in H^1(\mathcal{O})$, then $\mathcal{M} u \in H^1(\mathbb R^n)$ and
$$
\|u-\mathcal{M} u\|_{H_\rho(\mathcal{O})} \leq e_1 \varepsilon\|u\|_{H_{\varepsilon}^1(\mathcal{O})},
$$
where $e_1>0$ is a constant independent of $\varepsilon$ and $u$.
\end{lm}
To that end, we assume that all the functions $f_{\varepsilon},\varpi_{k}$ in \eqref{c36} satisfy the conditions $\mathbf(A1)$-$\mathbf(A3)$.
Furthermore, we assume that all $\varepsilon\in(0.1)$, $t,s\in\mathbb R$, $y^*\in\mathbb R^n$ and $\mu\in\mathcal P_2(L^2(\mathcal O))$,
\begin{align}\label{e4}
|f_{\varepsilon}(y,s,\mu)-f_0(y^*,s,\mu)|\le \varepsilon\kappa_1,\quad  \text{for all}\ s\in\mathbb R,
\end{align}
\begin{align}\label{e5}
\|\sigma_{k,\varepsilon}-\sigma_{k,0}\|_{L^2(\mathcal O)}\le \varrho_k\varepsilon,
\end{align}
 where $\kappa_1>0,$ and $\varrho_k,k\in\mathbb N$, with $\sum\limits^{\infty}_{k=1}\xi_k^2<\infty$.

We now write the process associated with \eqref{c8} as $U^{g,\varepsilon}$ and use $U^{g,0}$ for the process associated with \eqref{c9}.
The uniform measure attractors  of  $U^{g,\varepsilon}$ and $U^{g,0}$ are denoted by  $\mathcal A_{\varepsilon}$ and $\mathcal A_{0}$, respectively.

\begin{lm}\label{e7}
Suppose $\mathbf(A1)$-$\mathbf(A3)$, \eqref{c10} and \eqref{e4}-\eqref{e5} hold.   For every $\tau\in\mathbb R$, $T>0$ and $R>0$, and $g\in\mathcal H(g_0)$, we have
$$
\lim _{\varepsilon \rightarrow 0} \sup_{\mathbb E(\|\xi^{\varepsilon}\|^2_{H^1_{\varepsilon}(\mathcal{O})})\le R}\sup_{\tau\le t\le\tau+T}\mathbb E\left(\|u^{g,\varepsilon}(t,\tau,\xi^{\varepsilon})
-u^{g,0}(t,\tau,\mathcal M\xi^{\varepsilon})\|^2\right)=0.
$$
\begin{proof}
 Let $u^{\varepsilon}(t)=u^{g,\varepsilon}\left(t, \tau, \xi^{\varepsilon}\right)$, $u^{0}(t)=u^{g,0}\left(t, \tau, \mathcal{M} \xi^{\varepsilon}\right)$ and $v^{\varepsilon}(t)=u^{\varepsilon}(t)- u^{0}(t)$. By (3.24) and (3.25), we have, for all $t>\tau$,

$$
\begin{aligned}
& d v^{\varepsilon}(t)+\left(A_{\varepsilon} u^{\varepsilon}(t)-A_0 u^0(t)\right) d t+\lambda v^{\varepsilon}(t) d t+\left(f_{\varepsilon}\left(y,  u^{\varepsilon}(t),\mathcal L_{u^{\varepsilon}}\right)-f_0\left(y^*, u^0(t),\mathcal L_{u^{\varepsilon}}\right)\right) d t \\
& =
\sum_{k=1}^{\infty}\left(\left(\sigma_{k, \varepsilon}(y)-\sigma_{k, 0}\left(y^*\right)\right)+\left(\varpi_k\left(u^{\varepsilon}(t),\mathcal L({u^{\varepsilon}(t)})\right)-\varpi_k\left(u^0(t).\mathcal L({u^{0}(t)}\right)\right)\right) d W_k(t).
\end{aligned}
$$
By Ito's formula, we obtain for $t \geq \tau$,
\begin{align}\label{f2}
\begin{split}
& \left\|v^{\varepsilon}(t)\right\|_{H_\rho(\mathcal{O})}^2 \\
& =\left\|\xi^{\varepsilon}-\mathcal M\xi^{\varepsilon}\right\|_{H_\rho(\mathcal{O})}^2-2 \int_\tau^t\left(A_{\varepsilon} u^{\varepsilon}(s)-A_0 u^0(s), v^{\varepsilon}(s)\right)_{H_\rho(\mathcal{O})} d s\\
&-2\lambda\int^t_{\tau}\|v^{\varepsilon}(s)\|_{\mathcal H_{\rho}(\mathcal O)}ds-2\int^t_{\tau}\left(f_{\varepsilon}\left(\cdot,  u^{\varepsilon}(s),\mathcal L_{u^{\varepsilon}(s)}\right)-f_0\left(\cdot, u^0(s),\mathcal L_{u^{\varepsilon}(s)}\right),v^{\varepsilon}\right)_{H_\rho(\mathcal{O})}ds\\
&+\sum_{k=1}^{\infty} \int_\tau^t\left\|\sigma_{k, \varepsilon}-\sigma_{k, 0}+\kappa\varpi_k\left(u^{\varepsilon}(s),\mathcal L_{u^{\varepsilon}(s)}\right)-\kappa\varpi_k\left(u^0(s),\mathcal L_{u^{0}(s)}\right)\right\|_{H_\rho(\mathcal{O})}^2 d s \\
& +2 \sum_{k=1}^{\infty} \int_\tau^t\left(\sigma_{k, \varepsilon}-\sigma_{k, 0}, v^{\varepsilon}(s)\right)_{H_\rho(\mathcal{O})} d W_k(s) \\
& +2 \sum_{k=1}^{\infty} \int_\tau^t\left(\kappa\varpi_k\left(u^{\varepsilon}(s),\mathcal L_{u^{\varepsilon}(s)}\right)-\kappa\varpi_k\left(u^0(s),\mathcal L_{u^{0}(s)}\right), v^{\varepsilon}(s)\right)_{H_\rho(\mathcal{O})} d W_k(s).
\end{split}
\end{align}
For the second term on the right-hand of \eqref{f2}, we have
\begin{align}\label{f3}
\begin{split}
&-2 \int_\tau^t\left(A_{\varepsilon} u^{\varepsilon}(s)-A_0 u^0(s), v^{\varepsilon}(s)\right)_{H_\rho(\mathcal{O})} d s\\
&=-2\int_\tau^t\left(A_{\varepsilon} u^{\varepsilon}(s),v^{\varepsilon}(s)\right)_{H_\rho(\mathcal{O})} ds+2\int_\tau^t\left(A_{0} u^{0}(s),v^{\varepsilon}(s)\right)_{H_\rho(\mathcal{O})} ds\\
&=-2\int_\tau^t\left(A_{\varepsilon} u^{\varepsilon}(s),v^{\varepsilon}(s)\right)_{H_\rho(\mathcal{O})} ds+2\sum^{n}_{i=1}\int_\tau^t\left(u^0_{y_i}(s),v^{\varepsilon}_{y_i}(s)\right)_{H_\rho(\mathcal{O})}\\
&=-2\int_\tau^ta_{\varepsilon}\left(u^{\varepsilon}(s),v^{\varepsilon}(s)\right)ds+2\int_\tau^ta_{\varepsilon}\left(u^{0}(s),v^{\varepsilon}(s)\right)ds\\
&+2\sum^{n}_{i=1}\int^t_{\tau}(\frac{\rho_{y_i}}{\rho}u^0_{y_i}(s),y_{n+1}(u^{\varepsilon}_{y_{n+1}}(s)-u^0_{y_{n+1}}(s)))_{H_{\rho}(\mathcal{O})}ds\\
&=-2\int^t_{\tau}a_{\varepsilon}\left(v^{\varepsilon}(s),v^{\varepsilon}(s)\right)ds\\
&+2\sum^{n}_{i=1}\int^t_{\tau}(\frac{\rho_{y_i}}{\rho}u^0_{y_i}(s),y_{n+1}(u^{\varepsilon}_{y_{n+1}}(s)-u^0_{y_{n+1}}(s)))_{H_{\rho}(\mathcal{O})}ds.
\end{split}
\end{align}
For the fourth term on the right-hand of \eqref{f2}, by \eqref{c4}-\eqref{c14} and \eqref{e4} we have
\begin{align}
\begin{split}
&-2\int^t_{\tau}\int_{\mathcal O}\rho\left(f_{\varepsilon}\left(y,  u^{\varepsilon}(s),\mathcal L_{u^{\varepsilon}(s)}\right)-f_0\left(y^*, u^0(s),\mathcal L_{u^{0}(s)}\right)v^{\varepsilon}(s,y)\right)dyds\\
&=-2\int^t_{\tau}\int_{\mathcal O}\rho (f(y^*,\varepsilon g(y^*)y_{n+1}, u^{\varepsilon}(s,y),\mathcal L_{u^{\varepsilon}(s)})-f(y^*,\varepsilon g(y^*)y_{n+1}, u^{0}(s,y),\mathcal L_{u^{\varepsilon}(s)}))v^{\varepsilon}(s,y)dyds\\
&-2\int^t_{\tau}\int_{\mathcal O}\rho (f(y^*,\varepsilon g(y^*)y_{n+1}, u^{0}(s,y),\mathcal L_{u^{\varepsilon}(s)})-f(y^*,\varepsilon g(y^*)y_{n+1}, u^{0}(s,y),\mathcal L_{u^{0}(s)}))v^{\varepsilon}(s,y)dyds\\
&-2\int^t_{\tau}\int_{\mathcal O}\rho f(y^*,\varepsilon g(y^*)y_{n+1}, u^{0}(s,y),\mathcal L_{u^{0}(s)}))-f(y^*,0, u^{0}(s,y),\mathcal L_{u^{0}(s)}))v^{\varepsilon}(s,y)dyds\\
&\le2\int^t_{\tau}\int_{\mathcal O}\rho \phi_4(y^*)|v^{\varepsilon}(s,y)|^2dyds+2\int^t_{\tau}\int_{\mathcal O}\rho\phi_3(y^*)|v^{\varepsilon}(s,y)|\sqrt{\mathbb E(\|v^{\varepsilon}(s)\|^2)}dyds\\
&+\varepsilon\left(\kappa_1^2(t-\tau)+\int^t_{\tau}\|v^{\varepsilon}(s)\|^2_{H_{\rho}(\mathcal{O})} ds\right)\\
&\le2\int^t_{\tau}\|\phi_4\|_{L^{\infty}(\mathbb R^n)}\|v^{\varepsilon}(s)\|^2_{H_{\rho}(\mathcal{O})}ds\\
&+\int^t_{\tau}(\|\phi_3\|_{L^{\infty}(\mathbb R^n)}\|v^{\varepsilon}(s)\|^2_{H_{\rho}(\mathcal{O})}+\|\phi_3\|_{L^{1}(\mathbb R^n)}\mathbb E(\|v^{\varepsilon}(s)\|^2_{H_{\rho}(\mathcal{O})}))ds\\
&+\varepsilon\left(\kappa_1^2(t-\tau)+\int^t_{\tau}\|v^{\varepsilon}(s)\|^2_{H_{\rho}(\mathcal{O})} ds\right).
\end{split}
\end{align}
By \eqref{e5}, we have
\begin{align}
\sum^{\infty}_{k=1}\int^t_{\tau}\|\sigma_{k,\varepsilon}-\sigma_{k,0}\|^2_{H_\rho(\mathcal{O})}dt\le c_1\varepsilon^2(t-\tau).
\end{align}
It follows from \eqref{c7} we have
\begin{align}
&\|\kappa\left(\varpi_k\left(u^{\varepsilon}(s),\mathcal L_{u^{\varepsilon}(s)}\right)-\varpi_k\left(u^0(s),\mathcal L_{u^{0}(s)}\right)\right)\|_{H_\rho(\mathcal{O})}^2\nonumber\\
&\le2\|L_{\varpi}\|^2_{l^2}\int^t_{\tau}\left(\|\kappa\|^2_{L^{\infty}(\mathbb R^n)}\|v^{\varepsilon}(s)\|^2_{H_\rho(\mathcal{O})}+\|\kappa\|_{L^2(\mathbb R^n)}
\mathbb E(\|v^{\varepsilon}(s)\|^2_{H_\rho(\mathcal{O})})\right)ds.
\end{align}
Nexy, by \eqref{c19}, we obtain
\begin{align}\label{f4}
\begin{split}
&2\sum^{n}_{i=1}\int^t_{\tau}(\frac{\rho_{y_i}}{\rho}u^0_{y_i}(s),y_{n+1}(u^{\varepsilon}_{y_{n+1}}(s)-u^0_{y_{n+1}}(s))))_{H_{\rho}(\mathcal{O})}ds\\
&\le2\sum^{n}_{i=1}\int^t_{\tau}(\rho_{y_i}u^0_{y_i}(s),y_{n+1}(u^{\varepsilon}_{y_{n+1}}(s)-u^0_{y_{n+1}}(s)))_{L^2(\mathcal O)}ds\\
&\le c_2\varepsilon\int^t_{\tau}\left(\|u^0(s)\|^2_{H^1(\mathbb R^n)}\|u^{\varepsilon}(s)-u^0(s)\|_{H^1_{\varepsilon}(\mathcal O)}\right)ds\\
&\le c_3\varepsilon\int^t_{\tau}\|u^0(s)\|^2_{H^1(\mathbb R^n)}+\|u^{\varepsilon}(s)\|^2_{H^1_{\varepsilon}(\mathcal O)}ds.
\end{split}
\end{align}
Taking the expectation of \eqref{f2} and using \eqref{f3}-\eqref{f4}, we obtain for all $t\ge\tau$ and $\varepsilon\in(0,\varepsilon_0)$,
\begin{align}\label{f5}
\begin{split}
&\mathbb E\left(\left\|v^{\varepsilon}(t)\right\|_{H_\rho(\mathcal{O})}^2\right)\le
\mathbb E\left(\left\|\xi^{\varepsilon}-\mathcal M\xi^{\varepsilon}\right\|_{H_\rho(\mathcal{O})}^2\right)\\
&+c_3\varepsilon\int^t_{\tau}\mathbb E\left(\|u^0(s)\|^2_{H^1(\mathbb R^n)}+\|u^{\varepsilon}(s)\|^2_{H^1_{\varepsilon}(\mathcal O)}\right)ds\\
&+c_4\int^t_{\tau}\mathbb E\left(\|v^{\varepsilon}(s)\|^2_{H_{\rho}(\mathcal{O})}\right)ds+\varepsilon\kappa_1^2(t-\tau)
+c_1\varepsilon^2(t-\tau).
\end{split}
\end{align}
By \eqref{f5}, Lemma \ref{d38} and Lemma \ref{d39} we find that for every $T>0$, there exists $c_5=c_5(T,R)>0$ such that for all
$\tau\le T\le \tau+T$, $\varepsilon\in(0,\varepsilon_0)$ and $\xi^{\varepsilon}\in L^2(\Omega,\mathcal F_{\tau}, L^2(\mathcal O))$ with
$\mathbb E(\|\xi^{\varepsilon}\|^2_{H^1_{\varepsilon}(\mathcal O)})\le R$,
\begin{align}
\begin{split}
&\mathbb E\left(\left\|v^{\varepsilon}(t)\right\|_{H_\rho(\mathcal{O})}^2\right)\le
\mathbb E\left(\left\|\xi^{\varepsilon}-\mathcal M\xi^{\varepsilon}\right\|_{H_\rho(\mathcal{O})}^2\right)\\
&+c_4\int^t_{\tau}\mathbb E\left(\|v^{\varepsilon}(s)\|^2_{H_{\rho}(\mathcal{O})}\right)ds+\varepsilon\kappa_1^2T+c_5\varepsilon
+c_1\varepsilon^2T.
\end{split}
\end{align}
Then by Gronwall's inequality and Lemma \ref{f6}, we infer that for all $\tau\le t\le\tau+T$, $\varepsilon\in(0,\varepsilon_0)$
and $\xi^{\varepsilon}\in L^2(\Omega,\mathcal F_{\tau}, L^2(\mathcal O))$ with
$\mathbb E(\|\xi^{\varepsilon}\|^2_{H^1_{\varepsilon}(\mathcal O)})\le R$,
\begin{align}\label{f7}
\begin{split}
&\mathbb E\left(\left\|v^{\varepsilon}(t)\right\|_{H_\rho(\mathcal{O})}^2\right)\le
\left(\mathbb E\left(\left\|\xi^{\varepsilon}-\mathcal M\xi^{\varepsilon}\right\|_{H_\rho(\mathcal{O})}^2\right)+\varepsilon\kappa_1^2T+c_5\varepsilon
+c_1\varepsilon^2T\right)e^{c_4(t-\tau)}\\
&\le\left(c_6\mathbb E\left(\left\|\xi^{\varepsilon}\right\|_{H^1_{\varepsilon}(\mathcal{O})}^2\right)+\varepsilon\kappa_1^2T+c_5\varepsilon
+c_1\varepsilon^2T\right)e^{c_4T}\\
&\le\left(c_6R\varepsilon+\varepsilon\kappa_1^2T+c_5\varepsilon
+c_1\varepsilon^2T\right)e^{c_4T}.
\end{split}
\end{align}
It follows from \eqref{f7} that
\begin{align}\label{f8}
\lim _{\varepsilon \rightarrow 0} \sup_{\mathbb E(\|\xi^{\varepsilon}\|^2_{H^1_{\varepsilon}(\mathcal{O})})\le R}\sup_{\tau\le t\le\tau+T}\mathbb E\left(\|u^{g,\varepsilon}(t,\tau,\xi^{\varepsilon})
-u^{g,0}(t,\tau,\mathcal M\xi^{\varepsilon})\|^2_{L^2(\mathcal O)}\right)=0.
\end{align}
\end{proof}
\end{lm}

\begin{cor}\label{f9}
Suppose $\mathbf(A1)$-$\mathbf(A3)$, \eqref{c10} and \eqref{e4}-\eqref{e5} hold.   For every $\tau\in\mathbb R$, $t\ge\tau$ and $R>0$, and $g\in\mathcal H(g_0)$, we have
$$
\lim _{\varepsilon \rightarrow 0} \sup_{\mu^{\varepsilon}\in B^1_{\varepsilon}(R)}d_{\mathcal P(L^2(\mathcal O))}
\left(U^{g,\varepsilon}(t,\tau)\mu^{\varepsilon},\left(U^{g,0}(t,\tau)(\mu^{\varepsilon}\circ\mathcal M^{-1})\right)\circ\mathcal I^{-1}\right)=0,
$$
where $B^1_{\varepsilon}(R)=\{\mu\in\mathcal P_2(L^2(\mathcal O)):\int_{{H^1_{\varepsilon}(\mathcal{O})}}\|\xi\|^2_{{H^1_{\varepsilon}(\mathcal{O})}}\mu(d\xi)\le R\}$.
\begin{proof}
Note that for all $t\ge\tau$ we have
\begin{align}
\begin{split}
&\sup_{\mathbb E(\|\xi^{\varepsilon}\|^2_{H^1_{\varepsilon}(\mathcal{O})})\le R}\sup_{\substack{\varphi\in L_b(L^2(\mathcal O))\\ \|\varphi\|_L\le 1}}
|\mathbb E\left(\varphi(u^{g,\varepsilon}(t,\tau,\xi^{\varepsilon}))\right)- E\left(\varphi(u^{g,0}(t,\tau,\mathcal M\xi^{\varepsilon}))\right)|\\
&\le\sup_{\mathbb E(\|\xi^{\varepsilon}\|^2_{H^1_{\varepsilon}(\mathcal{O})})\le R}\sup_{\substack{\varphi\in L_b(L^2(\mathcal O))\\ \|\varphi\|_L\le 1}}
\mathbb E\left(|\varphi(u^{g,\varepsilon}(t,\tau,\xi^{\varepsilon}))- \varphi(u^{g,0}(t,\tau,\mathcal M\xi^{\varepsilon}))|\right)\\
&\le\sup_{\mathbb E(\|\xi^{\varepsilon}\|^2_{H^1_{\varepsilon}(\mathcal{O})})\le R}
\mathbb E\left(\|(u^{g,\varepsilon}(t,\tau,\xi^{\varepsilon}))- (u^{g,0}(t,\tau,\mathcal M\xi^{\varepsilon}))\|_{L^2(\mathcal O)}\right)\\
&\le\left(\sup_{\mathbb E(\|\xi^{\varepsilon}\|^2_{H^1_{\varepsilon}(\mathcal{O})})\le R}
\mathbb E\left(\|(u^{g,\varepsilon}(t,\tau,\xi^{\varepsilon}))- (u^{g,0}(t,\tau,\mathcal M\xi^{\varepsilon}))\|^2_{L^2(\mathcal O)}\right)\right)^{\frac{1}{2}},
\end{split}
\end{align}
which along with \eqref{f8} implies that for all $t\ge\tau$,
\begin{align}
\lim_{\varepsilon\rightarrow0}\sup_{\mathbb E(\|\xi^{\varepsilon}\|^2_{H^1_{\varepsilon}(\mathcal{O})})\le R}\sup_{\substack{\varphi\in L_b(L^2(\mathcal O))\\ \|\varphi\|_L\le 1}}|\mathbb E\left(\varphi(u^{g,\varepsilon}(t,\tau,\xi^{\varepsilon}))\right)- E\left(\varphi(u^{g,0}(t,\tau,\mathcal M\xi^{\varepsilon}))\right)|=0.
\end{align}
This completes the proof.
\end{proof}
\end{cor}

Next, we discuss the the upper semicontinuity of uniform measure attractors of \eqref{c8}.
\begin{thm}
Assume that $\mathbf(A1)$-$\mathbf(A3)$ and \eqref{e4}-\eqref{e5} hold. Then the uniform measure attractors $A_{\varepsilon}$ are upper semicontinuous
at $\varepsilon=0$,
\begin{align}\label{e6}
\lim\limits_{\varepsilon\rightarrow0}d_{\mathcal P(L^2(\mathcal O))}\left(\mathcal A_{\varepsilon},\mathcal A_0\circ \mathcal{I}^{-1}\right)=0.
\end{align}
\begin{proof}
By Lemma \ref{d10} we find that
\begin{align}\label{f1}
\int_{H^{1}_{\varepsilon}(\mathcal O)}\|\xi\|^2_{H^{1}_{\varepsilon}(\mathcal O)}\mu(d\xi)\le K_1\quad \text{for all}\ 0<\varepsilon<\varepsilon_0\ \text{and}\ \mu\in\mathcal A_{\varepsilon},
\end{align}
where $K_1>0$ is independent of $\varepsilon,g$. Let $K$ be the uniform absorbing set of $U^{g_{},\varepsilon}$ as given by \eqref{d9}, and denote by
 $K_0=\{\mu\circ\mathcal M^{-1}:\mu\in K\}$. Since $\mathcal A_0$ is the uniform measure attractor of $\{U^{g,0}\}_{g_{}\in\mathcal H(g_0)}$ in $\mathcal P_4(L^2(\mathbb R^n))$,
 given $\eta>0$, we infer that there exists $T=T(\eta)\ge1$ such that for any $t-\tau>T$ and $g\in\mathcal H(g_0)$,
\begin{align}\label{e8}
d_{\mathcal P_2(L^2(\mathbb R^n))}(U^{g,0}(t,\tau)K_0,\mathcal A_0)<\frac{1}{2}\eta.
\end{align}
On the other hand, by \eqref{f1} and Corollary \ref{f9}  we have
\begin{align}
\lim _{\varepsilon \rightarrow 0} \sup _{\nu_{\varepsilon} \in \mathcal A_{\varepsilon}} d_{\mathcal{P}\left(L^2(\mathcal{O})\right)}\left(U^{g_{},\varepsilon}(t, \tau) \nu_{\varepsilon}, \ \left(U^{g,0}(t,\tau)\left(\nu_{\varepsilon} \circ \mathcal{M}^{-1}\right)\right) \circ \mathcal{I}^{-1}\right)=0.
\end{align}
and hence there exists $\varepsilon_1\in(0,\varepsilon_0)$ such that for all $0<\varepsilon<\varepsilon_1$,
\begin{align}\label{e9}
\sup _{\nu_{\varepsilon} \in \mathcal A_{\varepsilon}} d_{\mathcal{P}_2\left(L^2(\mathcal{O})\right)}\left(U^{g_{},\varepsilon}(t, \tau) \nu_{\varepsilon}, \ \left(U^{g,0}(t,\tau)\left(\nu_{\varepsilon} \circ \mathcal{M}^{-1}\right)\right) \circ \mathcal{I}^{-1}\right)<\frac{1}{2}\eta
.
\end{align}
Given $\nu_{\varepsilon}\in\mathcal A_{\varepsilon}$, since $\mathcal A_{\varepsilon}\in K$, we know $\nu_{\varepsilon}\circ\mathcal M^{-1}\in K_0$, and thus
by \eqref{e8} we have
\begin{align}
\sup_{\nu_\varepsilon\in\mathcal A_{\varepsilon}}d_{\mathcal P(L^2(\mathbb R^n))}(U^{g,0}(t,\tau)\left(\nu_{\varepsilon} \circ \mathcal{M}^{-1}\right),\mathcal A_0)<\frac{1}{2}\eta.
\end{align}
which shows that
\begin{align}\label{e10}
\sup_{\nu_\varepsilon\in\mathcal A_{\varepsilon}}d_{\mathcal P(L^2(\mathcal O))}\left(\left(U^{g,0}(t,\tau)(\nu_{\varepsilon} \circ \mathcal{M}^{-1}\right)\right)\circ\mathcal I^{-1},\mathcal A_0\circ \mathcal I^{-1})<\frac{1}{2}\eta.
\end{align}
By \eqref{e9} and \eqref{e10} we have, for all $0<\varepsilon<\varepsilon_1$,
\begin{align}\label{e11}
\sup_{\nu_\varepsilon\in\mathcal A_{\varepsilon}}d_{\mathcal P(L^2(\mathcal O))}(U^{g_{},\varepsilon}(t,\tau)\nu_{\varepsilon},\mathcal A_0\circ \mathcal I^{-1})<\eta.
\end{align}
By the uniformly quasi-invariance of $A_{\varepsilon}$, we see that for any $\mu_{\varepsilon}\in\mathcal A_{\varepsilon}$, there exists $\nu_{\varepsilon}\in\mathcal A_{\varepsilon}$ and $g\in \mathcal H(g_0)$ such that
\begin{align}\label{e12}
\mu_{\varepsilon}=U^{g_{},\varepsilon}(t,\tau)\nu_{\varepsilon}.
\end{align}
By \eqref{e11} and \eqref{e12} we obtain, for all $0<\varepsilon<\varepsilon_1$,
$$\sup_{\mu_\varepsilon\in\mathcal A_{\varepsilon}}d_{\mathcal P(L^2(\mathcal O))}(\mu_{\varepsilon},\mathcal A_0\circ \mathcal I^{-1})<\eta,$$
which indicates that for all $0<\varepsilon<\varepsilon_1$,
$$d_{\mathcal P(L^2(\mathcal O))}(\mathcal A_{\varepsilon},\mathcal A_0\circ \mathcal I^{-1})<\eta,$$
as desired.
\end{proof}
\end{thm}

\bibliographystyle{plain}
{\footnotesize


\end{document}